\documentclass[12pt,reqno]{article}

\usepackage[usenames]{color}
\usepackage{amssymb}
\usepackage{graphicx}
\usepackage{amscd}

\usepackage[colorlinks=true,
linkcolor=webgreen,
filecolor=webbrown,
citecolor=webgreen]{hyperref}

\definecolor{webgreen}{rgb}{0,.5,0}
\definecolor{webbrown}{rgb}{.6,0,0}

\usepackage{color}
\usepackage{float}

\usepackage{graphics,amsmath,amssymb}
\usepackage{amsthm}
\usepackage{amsfonts}
\usepackage{latexsym}

\newcommand{\seqnum}[1]{\href{http://oeis.org/#1}{\underline{#1}}}

\begin{document}

\theoremstyle{plain}
\newtheorem{theorem}{Theorem}
\newtheorem{corollary}[theorem]{Corollary}
\newtheorem{lemma}[theorem]{Lemma}
\newtheorem{proposition}[theorem]{Proposition}
\newtheorem{obs}[theorem]{Observation}

\theoremstyle{definition}
\newtheorem{definition}[theorem]{Definition}
\newtheorem{example}[theorem]{Example}
\newtheorem{conjecture}[theorem]{Conjecture}
\newtheorem{question}[theorem]{Question}

\theoremstyle{remark}
\newtheorem{remark}[theorem]{Remark}

\begin{center}
\vskip 1cm
{\Large\bf An alternate form of Merino-Mi\v{c}ka-M\"utze's solution to a Knuth's combinatorial generation problem}%

\vskip 1cm
\large
Italo J. Dejter\\
University of Puerto Rico\\
Rio Piedras, PR 00936-8377\\
\href{mailto:italo.dejter@gmail.com}{\tt italo.dejter@gmail.com} \\
\end{center}
 
\begin{abstract}\noindent 
A modification of Merino-Mi\v{c}ka-M\"utze's solution to a combinatorial generation problem of Knuth is proposed in this survey. The resulting alternate form to such solution is compatible with a reinterpretation by the author of a proof of existence of Hamilton cycles in the middle-levels graphs. Such reinterpretation is given in terms of a dihedral quotient graph associated to each middle-levels graph. The vertices of such quotient graph represent Dyck words and their associated ordered trees. Those Dyck words are linearly ordered via a rooted tree that covers all their tight, or irreducible, forms, offering an universal reference point of view to express and integrate the periodic paths, or blocks, whose concatenation leads to Hamilton cycles resulting from the said solution.
\end{abstract}

\section{Introduction}\label{s1} In this survey,
an alternate form of Merino-Mi\v{c}ka-M\"utze's solution \cite{mmm} to a combinatorial generation problem of Knuth \cite{Knuth} is given that is compatible with a reinterpretation \cite{D1} of the proof by M\"utze~\cite{M} and Gregor, M\"utze and Nummenpalo~\cite{gmn} 
of existence 
of Hamilton cycles in the {\it middle-levels graph} $M_n$ \cite{KT} ($0<n\in\mathbb{Z}$). This graph $M_n$ is induced by the vertices of the $(2n+1)$-cube that represent the $n$- and $(n+1)$-subsets of the set $\{0,\ldots,2n\}$ and is further considered from Section~\ref{mlgs} on. 

The reinterpretation in question is given in terms of the associated dihedral quotient graph $N_n$ of $M_n$ whose vertices (here called {\it necklaces}, as in~\cite{mmm}, see Definition~\ref{defneck}, below) represent the Dyck words of length $2n$ defined in the following alternate way to that of~\cite{mmm}:

\begin{definition} The {\it deficiency} (not the {\it excess}, as in~\cite{mmm}) of a bitstring $x$ is its number of 0's minus its number of 1's. If $x$ has deficiency 0 and every prefix has negative deficiency, then we say that $x$ is a {\it Dyck word}.
Let $D_n$ be the set of Dyck words of length $2n$. Let $D = \cup_{n\ge 0} D_n$.
\end{definition}

Such words can be linearly ordered via a {\it castling}, or {\it lexical}, procedure \cite[Section 3]{D1} (where the term {\it lexical} appears in \cite{D1} in relation to Kierstead and Trotter~\cite{KT}).

The castling procedure uses {\it restricted growth strings} or RGSs of length $n$, namely the $n$-strings $a_1a_2\cdots a_n$ such that $a_1=0$ and $a_k\le a_{k-1}+1$, for $1<k\le n$; see \cite[page 325]{Arndt}-\cite[page 224 (u)]{Stanley}, that establish that their number is the Catalan number $C_n=\frac{1}{n+1}{2n\choose n}$ \cite[\seqnum{A100108}]{oeis}. By eliminating the initial zero of each such RGS, we say that the resulting $(n-1)$-string is an $n$-{\it germ} $a_2\cdots a_n=b_{n-1}\cdots b_0$, better denoted $b_{n-1}\cdots b_0$ in what follows.

Those $n$-germs form an ordered tree $\mathcal T^n$ \cite[Theorem 1]{D1} whose root is the null $n$-germ $\alpha=0^{n-1}$ and such that each non-null $n$-germ $\alpha=a_{n-1}\cdots a_1$ with rightmost nonzero entry $a_i$ ($1\le i<n$) has as parent in $\mathcal T_n$ the $n$-germ $\beta=b_{n-1}\cdots b_1$ such that $b_i=a_i-1$ and $a_j=b_j$, for $j\ne i$.

To each non-null $n$-germ $\alpha$ corresponds its {\it tight RGS}, or {\it irreducible RGS}, obtained by removing its largest null prefix.
The ordered trees $\mathcal T^k$ form a chain ${\mathcal T^1}\subset{\mathcal T^2}\subset\cdots\subset{\mathcal T^n}\subset\cdots$ etc. of ordered trees obtained
by identifying each $n$-germ $\alpha$ with the $(n+1)$-germ obtained from $\alpha$ by prefixing a zero to it. The limit $\cup_{i=1}^{\infty}{\mathcal T^i}$ of this chain, denoted $\mathcal T$, will be considered with each vertex $v$ written as the tight RGS $\alpha$ associated to the $n$-germs $\alpha,0\alpha,00\alpha,\ldots$ that $v$ represents, so we have 
$$V({\mathcal T^1})=\{\emptyset\}, V({\mathcal T^2})=\{0,1\}, V({\mathcal T^2})=\{00,01,10,11,12\},\mbox{ etc.,}$$ that, as subsets of $\mathcal T$, will be expressed as
$$V({\mathcal T^1})=\{\emptyset\}, V({\mathcal T^2})=\{\emptyset,1\}, V({\mathcal T^2})=\{\emptyset,1,10,11,12\}\mbox{ etc.}$$
The tree $\mathcal T$ is partially represented on the lower-right of Table~\ref{table} by its subtree ${\mathcal T}^4$.

The castling procedure is realized by an inductively defined bijection $F$ from the set of $n$-germs onto $D_n$. This yields a natural way of designating the participating necklaces.

We ``personalize" each necklace of length $2n+1$ by taking a representation $f_0f_1\cdots f_{2n}$ of it such that $f_0=0$ and $f_1\cdots f_{2n}\in D_n$, and then applying the following procedure:
\begin{eqnarray}\label{ifni}\begin{array}{l}
g_0:=0;\\
\mbox{begin}\\
\hspace*{3mm}\mbox{ for }k=1\mbox{ to }2k\mbox{ do}\\
\hspace*{3mm}\mbox{if }f_{k+1}=0\mbox{ then }g_{k+1}=g_k+1\mbox{ else }g_{k+1}=g_k-1;\\
 k:=k+1\\
 \mbox{end.}
 \end{array}\end{eqnarray} 
 \noindent The procedure (\ref{ifni}) replaces each bit $f_k$ of the necklace $f_0f_1\cdots f_{2n}$ with a corresponding integer $g_k$ ($0\le g_k\le n$). We say that the resulting $(2n+1)$-tuple $g_0g_1\cdots g_{2n}$ is the {\it neck}, more specifically the $n$-{\it neck}, associated to the necklace $f_0f_1\cdots f_{2n}$, (with the $2n$-string $f_1\cdots f_{2n}$ denoted as an $n$-{\it nest} in~\cite{universal,castling}). It happens that the necklace $f_0f_1\cdots f_{2n}$ has each bit $f_k$ replaced by the height  $g_k$ it determines in its {\it Dyck path} \cite[Subsection 3.2]{D2}, so that each 0-bit is replaced either by the null height at the starting point of the Dyck path or the first appearance of each non-null height, and each 1-bit is replaced by the second appearance of that non-null height. See Remark~\ref{remark} and Figure~\ref{gris}, below. 
This ``personalized" representation of the necklaces, that we call necks here, are shown as an example of Theorem~\cite[Theorem 2]{D1} on the left of Table~\ref{table}, below, generating $\mathcal T^4$ via the procedure (\ref{ifni}), where each row, except the top one, shows an inductively obtained parent correspondence  
$F:n$-germ$^{parent}\mapsto n$-neck$^{parent}$ assigned to its child correspondence $F:n$-germ$^{child}\mapsto n$-neck$^{child}$.

In Table~\ref{table}, from each parent $n$-germ to its child $n$-germ, the  suffix in bold, shows a unit increase at the rightmost nonzero entry of the child's $n$-germ. Such suffix in bold has the same length as both the prefixes and suffixes in bold of both the parent neck and child neck. Both are to be kept unchanged, while a central-string castling takes place on the non-bold central strings. In fact, the one on the parent side is subdivided into a left Italic substring and a right Roman substring, while the other one, on the child side (i.e., inside $n$-neck$^{child}$), permutes the order of the two substrings, justifying denoting this procedure as ``castling", like in chess. 
The Roman substring of a central string starts at the first appearance of the integer $\ell+1$ subsequent to the integer $\ell$ starting the Italic side.

The upper-right of Table~\ref{table} represents how $n$-germs and $n$-necks fit for each index $n$ into the case of the subsequent index $n+1$, where each zero which is added as a prefix (between auxiliary parentheses) of an RGS results in the corresponding neck into an inserted substring (also between auxiliary parentheses) \cite{universal}. As a result,
the ordered tree $\mathcal T$ of RGSs yields a bijection into an ordered tree $\mathcal T'$ covering all {\it tight} necks (namely, where all parenthesized substrings as mentioned above are eliminated).

This section is completed by providing some historical motivation and statement of the main results of \cite{mmm} as Theorems~\ref{the} and~\ref{extra}. In the rest of the paper, we adapt the discourse of~\cite{mmm} to our alternate form. In Section~\ref{mlgs} we provide additional information about the middle-levels graphs and necklace graphs. In Section~\ref{flip}, basic flip sequences are presented. In Section~\ref{per}, periodic paths and gluing pairs are introduced. In Section~\ref{1st}, an initial attempt at proving theorem~\ref{the} is given. In Sections~\ref{sketch1} and~\ref{sketch2}, sketches of the proofs of Theorems~\ref{the} and~\ref{extra} are given, while Section~\ref{redef} provides information about involved efficient computations.

\begin{table}[htp]
$$\begin{array}{||r|rcl|rcl|rcl||}\hline\hline
0&&F&&00{\bf 0}&\rightarrow&012344321&(0)(0)(0)\emptyset&\rightarrow&01(2(3(44)3)2)1\\
1&00{\bf 0}&\rightarrow&{\bf 0}{\it 1}234432{\bf 1}&00{\bf 1}&\rightarrow&{\bf 0}234432{\it 1}{\bf 1}&(0)(0)1&\rightarrow&02(3(44)3)211\\
2&0{\bf00}&\rightarrow&{\bf 01}{\it 2}3443{\bf 21}&0{\bf 10}&\rightarrow&{\bf 01}3443{\it 2}{\bf 21}&(0)10&\rightarrow&013(44)3221\\
3&01{\bf 0}&\rightarrow&{\bf 0}{\it 13443}22{\bf 1}&01{\bf 1}&\rightarrow&{\bf 0}22{\it 13443}{\bf 1}&(0)11&\rightarrow&02213(44)31\\
4&01{\bf 1}&\rightarrow&{\bf 0}{\it 221}3443{\bf 1}&01{\bf 2}&\rightarrow&{\bf 0}3443{\it 221}{\bf 1}&(0)12&\rightarrow&03(44)32211\\
5&{\bf 000}&\rightarrow&{\bf 012}{\it 3}44{\bf 321}&{\bf 100}&\rightarrow&{\bf 012}44{\it 3}{\bf 321}&====&==&========\\
6&10{\bf 0}&\rightarrow&{\bf 0}{\it 1}244332{\bf 1}&10{\bf 1}&\rightarrow&{\bf 0}244332{\it 1}{\bf 1}&&_\swarrow\emptyset_{\downarrow\searrow}&\\
7&1{\bf 00}&\rightarrow&{\bf 01}{\it 244}33{\bf 21}&1{\bf 10}&\rightarrow&{\bf 01}33{\it 244}{\bf 21}&1.&10_\downarrow&100_{\downarrow\searrow}\\
8&11{\bf 0}&\rightarrow&{\bf 0}{\it 133}2442{\bf 1}&11{\bf 1}&\rightarrow&{\bf 0}2442{\it 133}{\bf 1}&&11_\downarrow&101.\;110_{\downarrow\searrow}\\
9&11{\bf 1}&\rightarrow&{\bf 0}{\it 2442133}{\bf 1}&112&\rightarrow&033{\it 24421}{\bf 1}&&12.&\hspace*{8mm}111_\downarrow\;120_\downarrow\\
10&1{\bf 10}&\rightarrow&{\bf 01}33{\it 244}{\bf 21}&1{\bf 20}&\rightarrow&{\bf 01}44{\it 332}{\bf 21}&&&\hspace*{8mm}112.\;121_\downarrow\\
11&12{\bf 0}&\rightarrow&{\bf 0}1443322{\bf 1}&12{\bf 1}&\rightarrow&{\bf 0}22{\it 14433}{\bf 1}&&&\hspace*{16mm}122_\downarrow\\
12&12{\bf 1}&\rightarrow&{\bf 0}{\it 22144}33{\bf 1}&12{\bf 2}&\rightarrow&{\bf 0}33{\it 22144}{\bf 1}&&&\hspace*{16mm}123.\\
13&12{\bf 2}&\rightarrow&{\bf 0}{\it 33221}44{\bf 1}&12{\bf 3}&\rightarrow&{\bf 0}44{\it 33221}{\bf 1}&&&\\\hline\hline
\end{array}$$
\caption{Castling for $n=4$, fitting cases $n=2,3,4$ and ordered subtree $\mathcal T^4$ of $\mathcal T$.}
\label{table}
\end{table}

\noindent{\bf Combinatorial generation:} 
An expressed objective of \cite{mmm} is to generate all $(k,\ell)${\it -com\-bi\-na\-tions}, i.e. all ways of choosing a subset $S$ of a fixed size $k$ from the set $[n]:=\{1,\ldots,n\}$, with $n=k+\ell$.
Each such subset $S$ is encoded by a bitstring of length $n$ with exactly $k$ many 1's, where the $i$-th bit is 1 if and only if the element $i$ is contained in $S$.
 
Buck and Wiedemann conjectured in \cite{BW} that all $(n+1,n+1)$-combinations are generated by {\it star transpositions}, for every $n\ge 1$, i.e. in each step the element 1 either enters or leaves the set. The corresponding {\it flip sequence} $\alpha$ records the position of the bit swapped with the first bit in each step, where positions are indexed in $[2n+1]$ and $\alpha$ has length $N:={2n+2\choose n+1}.$ Buck-Wiedemann's conjecture was independently raised by Havel \cite{Havel} and became known as the {\it middle-levels conjecture}, name coming from an equivalent formulation of the problem, which asks for a Hamilton cycle in the middle-levels graphs, recalled below in Section \ref{mlgs}.\\
 
 \noindent{\bf Knuth's conjecture:} In \cite[Problem 56, Section 7.2.1.3]{Knuth}, Knuth conjectured that there is a star transposition that orders the $(n+1,n+1)$-combinations, for every $n\ge 1$, such that the flip sequence $\alpha$ has a block structure $\alpha=(\alpha_0,\alpha_1,\ldots,\alpha_{2n})$, where each block $\alpha_i$  has length $\frac{N}{2n+1}=\frac{1}{2n+1}{2(n+1)\choose n+1}$ and is obtained from $\alpha_0$ by element-wise addition  of $i$ mod $2n+1$, where $i\in[2n]$. As the entries of $\alpha$ are from $[2n+1]$, the numbers $1,\ldots,2n+1$ are chosen as addition residue-class representatives, rather than the usual $0,\ldots,2n$. Note that $\frac{N}{2n+1}=2C_n$, where $C_n=\frac{1}{n+1}{2n\choose n}$ is the n-th Catalan number. Then, \cite{mmm} proves the following.
 
 \begin{theorem}\label{the}\cite[Theorem 1]{mmm}
 For any $n\ge 1$ and $1\le s\le 2n$ that is coprime to $2n+1$, there is a star-transposition ordering of all $(n+1,n+1)$-combinations such that the corresponding flip sequence is of the form $\alpha=(\alpha_0,\alpha_1,\ldots,\alpha_{2n})$ with each block $\alpha_i$ obtained from $\alpha_0$ by element-wise addition of $i.s$ modulo $2n+1$, where $i\in[2n]$.
 \end{theorem}
 
 \begin{remark}\label{pp} By omitting the first entry of every $(n+1,n+1)$-combination, the $(n+1,n+1)$-combinations are transformed bijectively into the vertices of the middle-levels graphs $M_n$, so Theorem~\ref{the} can be rephrased in terms of Hamilton cycles of $M_n$, each cycle formed as a concatenation of $2n+1$ copies of a periodic path representing a block as in the statement of the theorem, presented below in terms of the approach of \cite{D1}, with lemmas and propositions leading to the proof of Theorem~\ref{the} stated in parallel to those of \cite{mmm}.\end{remark}
 
 The proof of Theorem~\ref{the} in \cite{mmm} is constructive and translates into an algorithm that generates all $(n+1,n+1)$-combinations by star transpositions efficiently, stated in \cite{mmm} as follows.
 \begin{theorem}\label{extra}\cite[Theorem 2]{mmm} There is an algorithm that computes for any $n\ge 1$ and $1\le s\le 2n$ that is coprime to $2n+1$, a star transposition ordering of all $(n+1,n+1)$-combinations as in Theorem~\ref{the}, with running time $O(n)$ for each generated combination, using $O(n)$ memory.\end{theorem} 

Since $\mathcal T^n$ is an ordered tree, then its vertex set $V({\mathcal T_n})$ inherits a natural linear order $\mathcal L^n$. Thus, we have a chain of linear orders
$\mathcal L^1\subset\mathcal L^2\subset\cdots\subset\mathcal L^n\subset\cdots$ etc. The limit $\cup_{i=1}^\infty\mathcal L^i$ of such a chain, denoted $\mathcal L$ and consisting of tight RGSs, induces via the castling operation a linear order $\mathcal L'$ of the tight necks,
offering an universal reference point of view \cite{universal}--\cite{castling} to express and integrate the periodic paths or blocks whose concatenation leads to Hamilton cycles resulting from Remark~\ref{pp}.
 
\begin{question} Can the Hamilton cycles obtained in the middle-levels graphs, or the periodic paths found for Knuth problem, according to Remark~\ref{pp}, be set in terms of the tight RGSs in the linear order $\mathcal L$, or their numerical designations via the linear order $\mathcal L'$?
 \end{question}
 
\begin{question} Is it possible to find a periodic path $P_n$ in each $M_n$, for $n>1$, so that when considering the composing vertices of such paths as elements of $V({\mathcal T'})$ it happens that $P_n$ is an initial subpath (prefix) of $P_{n+1}$, $\forall n>1$?
\end{question}
 
\section{Middle-levels graphs and necklace graphs}\label{mlgs}
Let $0<n\in\mathbb{Z}$. 
Let $A_n$ (resp. $B_n$) be the set of bitstrings of length $2n+1$ and weight $n$ (resp. $n+1$).
The middle-levels graph $M_n$ can be set as the graph whose vertex set is $V(M_n)=A_n\cup B_n$ and whose adjacency is given by a single flip.
The positions of the bitstrings in $V(M_n)$ are denoted $1,2,...,2n+1$ (mod $2n+1$), where $2n+1$ takes the place of additive identity 0.
Let $\sigma^i(x)$ denote the cyclic right-rotation by $i$ positions (\cite{mmm} uses left-rotation, while our approach here is compatible with the treatment of \cite{D1}).
The {\it necklace} $\langle x\rangle$ of $x$ is defined to be $\{\sigma^i(x);i\ge 0\}$.
For example,
if $x=11000\in A_2$ then
$\langle x\rangle = \{11000, 01100, 00110, 00011, 10001\}$.\\

\begin{definition}\label{defneck} Define the {\it necklace graph} $N_n$ to have as vertex set all necklaces $\langle x\rangle$, ($x\in V(M_n)$), with an edge between $\langle x\rangle$ and $\langle y\rangle$ if and only if  $x$ and $y$ differ in a single bit. 
$N_n$ is quotient graph of $M_n$ under the equivalence relation given by cyclically rotating bitstrings.
There may be, for each $\langle x\rangle$, two distinct bits in $x$ that reach the same $\langle y\rangle$.
But $N_n$ is to be considered as a simple graph, so in $N_n$ not all vertices have the same degree.
$N_n$ has less vertices than  $M_n$ by a factor of $2n+1$.\end{definition}

\begin{definition}\label{periodic} To obtain a flip sequence for a Hamilton cycle in $M_n$, we say that 
a path $P=\{x_1,...,x_k\}$ in $M_n$ is {\it periodic} if flipping a single bit in $x_k$ yields a vertex $x_{k+1}$ that satisfies $\langle x_{k+1}\rangle = \langle x_1\rangle$.\end{definition}

\noindent{\bf Operations on sequences $x=(x_1,\ldots,x_n)$:}
of integers: $x+a:=(x_1+a,...,x_k+a)$, ($a\in\mathbb{Z}$) and $|x|=$ length of $x$; of bitstrings: $\langle x\rangle=(\langle x_1\rangle,\ldots,\langle x_k\rangle)\mbox{ and } \sigma^i(x)=(\sigma^i(x_1),\ldots,\sigma^i(x_k)).$\\

\begin{remark}\noindent{\bf Rooted trees:}\label{RT} Differing from \cite{mmm}, all {\it rooted trees} treated here have a specific right-to-left ordering for the children of each vertex. Every Dyck word $x\in D_n$ can be interpreted as one such rooted tree on $n$ edges, as follows, adapting the viewpoint of  \cite{mmm} to the setting of \cite{D1}, where
$\epsilon$ stands for the empty bitstring:
If $x=\epsilon$, then $x$ is associated to the tree formed by an isolated  root;
else, $x=u0v1$, where $u,v\in D$. The trees $R, L$ corresponding to $v, u$, respectively, have the tree corresponding to $x$ with $R$ rooted at the rightmost child of the root, and the edges from the root to all other children except that rightmost child, together with their subtrees, forming the tree $L$.
This yields a bijection from $D_n$ onto the rooted trees with $n$ edges.\end{remark}

\begin{figure}[htp]
\includegraphics[scale=0.7]{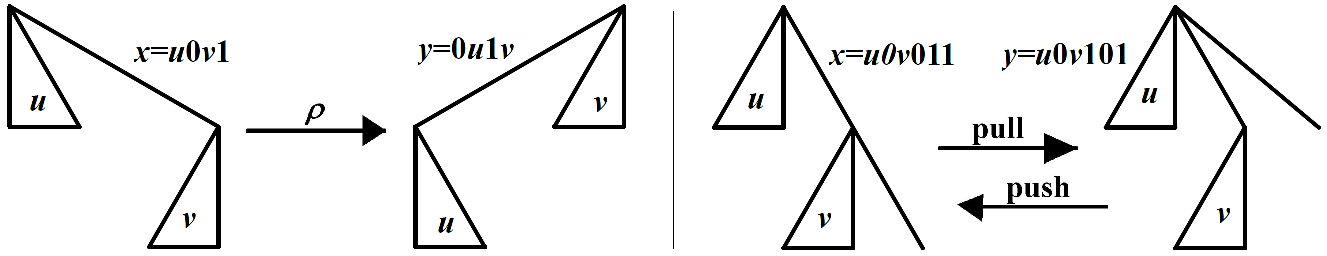}
\caption{Ordered-tree modifications}
\label{unalinea}
\end{figure}

\begin{remark}\label{rho}
{\bf Rooted-tree rotations:} Given a rooted tree $x\ne\epsilon$, let $\rho(x)$ denote the tree obtained by rotating $x$ to the left (in contrast to \cite{mmm}, that rotates it to the right), which corresponds to designating the rightmost child (not leftmost as in \cite{mmm}) of the root of $x$ as the root of $\rho(x)$. 
In terms of bitstrings, if $x = u0v1$, with $u,v\in D$, then $\rho(x) = 0u1v$. See the left half of Figure~\ref{unalinea}, which resembles, but differs reflectively from \cite[Figure 7]{mmm}.
\end{remark}

\begin{definition} A {\it plane tree} is a tree embedded in the plane with a specified clockwise cyclic ordering for the neighbors of each vertex, (not counterclockwise, shortened as ccw, \cite{mmm}).\end{definition} 

For $n\ge 1$, let $PT_n$ be the set of all plane trees with $n$ vertices. 
For any rooted tree $x$, let $[x]$ denote the set of all rooted trees obtained from $x$ by rotation (that is 
$[x] = \{\rho^i(x) ; i\ge 0\}$), 
to be interpreted as the plane tree underlying $x$, obtained by "forgetting" the root. 

Let $\lambda(x) = |[x]|$. For $T=[x]\in PT_n$, define 
$\lambda(T) = \lambda(x).$ Note that 
$$\lambda(x) = \min\{i\ge 1 ; \rho^i(x) = x\},$$ 
the choice of representative of $[x]$ in defining $\lambda(T)$ being irrelevant, for $\lambda(T)$ is well defined. Examples for $\lambda=4,8,2,3$ are given in Figure~\ref{gris}, below, in the notation of \cite{D1}, meaning that each 0-bit (resp. 1-bit) is represented by the first (resp. second) appearance of each integer, counting appearances rotationally from the red 0 and in the direction indicated by ``$>$" or by ``$<$". See also Remark~\ref{remark}.\\

\noindent{\bf $\bullet$-Subtrees:} Let $T\in PT_n$. Let $(a,b)\in E(T)$. Let $T^{(a,b)}$ be $T$ seen as a rooted tree with root $a$ and rightmost child $b$ (not leftmost as in \cite{mmm}). Let $T^{(a,b)-}$ be obtained from $T^{(a,b)}$ by removing all its children and their subtrees except for $b$ and its descendants.
 Given $a\in V(T)$, and all neighbors $b_i$ of $a$, ($i\in[k]$), let the trees $t_i=T^{(a,b_i)_-}$ be called the $a${\it -subtrees} of $T$. Then, $T=[(t_1,\ldots,t_k)]$, where $(t_1,\ldots,t_k)$ is the rooted tree obtained by gluing $t_1,\ldots,t_k$ at their roots from right to left (in this order, reversed to that of  \cite{mmm}):
 In terms of bitstrings, $(t_1,\ldots,t_k)$ is obtained by concatenating the bitstring representations of $t_1,\ldots,t_k$.\\

\noindent{\bf Centroids:} Given a (rooted or plane) tree $T$, the {\it potential} $\phi(a)$ of a vertex $a$ of $T$ is the sum of the distances from $a$ to $V(T)$. The {\it potential} $\phi(T)$ of $T$ is $\phi(T)=$min$\{\phi(a);a\in V(T)\}$. A {\it centroid} of $T$ is an $a\in V(T)$ with $\phi(a)=\phi(T)$. Merino, Mi\v{c}ka and M\"utze \cite{mmm}: {\bf(i)} mention that a centroid of $T$ is a vertex whose removal splits $T$ into subtrees with at most $\frac{|V(T)|}{2}$ vertices each;  {\bf(ii)}
prove in \cite[Lemma 3]{mmm} that $T$ has either one centroid or two adjacent centroids, and that if $|E(T))|$ is even, then $T$ has just one centroid.  

\begin{lemma}\label{ws}\cite[Lemma 4]{mmm} Let $T\in PT_n$ with $n\ge 1$ edges. Then, $\lambda(T)|2n$. If $T$ has a unique centroid, then $\lambda(T)$ is even; else, $\lambda(T) = 2n$ if $n$ is even, and $\lambda(T)\in\{n,2n\}$ if $n$ is odd. For $n\ge 4$ and any even divisor $k$ of $2n$, or for $k=n$, there is $T\in PT_n$ with $\lambda(T) = k$. \end{lemma}  

\noindent{\bf Relation of middle levels to Dyck words:} Our objective is to define as in \cite{mmm} basic flip sequences that visit every necklace exactly once in order to obtain a 2-factor (or cycle factor \cite{mmm}) in $N_n$, namely a collection of disjoint cycles that visits every vertex of $N_n$ exactly once.

\begin{lemma}\label{5}
Let $n\ge 1$. 
For any $x\in A_n$, there is a unique integer $\ell = \ell(x)$ with $0\le\ell\le 2n$ such that the last $2n$ bits of $\sigma^\ell(x)$ form a Dyck word.
For any $y\in B_n$, there is a unique integer $\ell = \ell(y)$ with $0\le\ell\le 2n$ such that the first $2n$ bits of $\sigma^\ell(y)$ form a Dyck word.
{\rm(Modified from \cite[Lemma 5]{mmm} that refers in turn to \cite[Problem 7]{Bol})}.
\end{lemma}

\noindent{\bf Dyck words of $A_n$ in the notation of Lemma~\ref{5}:} $\forall x\in A_n$, let $t(x)\in D_n$ denote the last $2n$ bits of $\sigma^\ell(x)$, where $\ell =\ell(x)$, i.e. $\sigma^\ell(x) = 0t(x)$, and $\forall y\in B_n$, let $t(y)\in D_n$ denote the first $2n$ bits of $\sigma^\ell(y)$, where $\ell =\ell(y)$, i.e. $\sigma^\ell(y) = t(y)1$.
Then, by Lemma~\ref{5}, every $x\in A_n$ (resp. $y\in B_n$) can be identified uniquely with the pair $(t(x),\ell(x))$ (resp. $(t(y),\ell(y))$).

\section{Basic flip sequences}\label{flip} 
A bijection $f$ on $V(M_n)$ is introduced
that yields a basic flip sequence visiting every necklace exactly once, however in a different fashion to that of \cite{mmm} but akin to the treatment of \cite{D1}.

Let $x\in A_n$ with $\ell(x) = 0$, i.e. 
$x = 0(u0v1) = 0t(x)$, where $u,v\in D$.
Define
 $y:=f(x) = (0u1v)1 = \rho(t(x))1\in B_n$,
where $\ell(y) = 1$ and $\rho$ is as in Remark~\ref{rho}.
Furhter, define 
$f(y) = f(f(x)) = (0u1v)0 = \rho(t(x))0\in A_n$, where $\ell(f(y)) =0$.
Extend these definitions of $f$ for all $x\in V(M_n)$ via 
$f(x) := \sigma^{-\ell}(f(\sigma^\ell))$, where $\ell := \ell(x)$.
Then $f$ is invertible and $t(f(f(x)))=t(f(x))=\rho(t(x))$. In our alternate situation (differing from that of~\cite{mmm}), $\ell(x)=1$, $\ell(f(x))=0$ and $\ell(f(f(x)))=0$. Then, for all $x\in A_n$, $\ell(f(x))=\ell(x)-1$ and $\ell(f(f(x)))=\ell(x)-1$.

\begin{definition}\label{p} For any $x\in V(M_n)$, let $\kappa(x)=\min\{i>0;\langle f^i(x)\rangle=\langle x\rangle\}$, be the {\it period} of $f$ at $x$, namely the number of  times $f$ must be applied before returning to the same necklace of $x$. 
For any $x\in V(M_n)$, let $P(x):=(x,f(x),f^2(x),\ldots,f^{\kappa(x)-1}(x))$ be the {\it periodic path} of $f$ at $x$ (Definition~\ref{periodic}) in $M_n$, namely a path of period $\kappa(x)$. Therefore, $\langle P(x)\rangle$ is a cycle in $N_n$.\end{definition}

\begin{figure}[htp]
\includegraphics[scale=0.7]{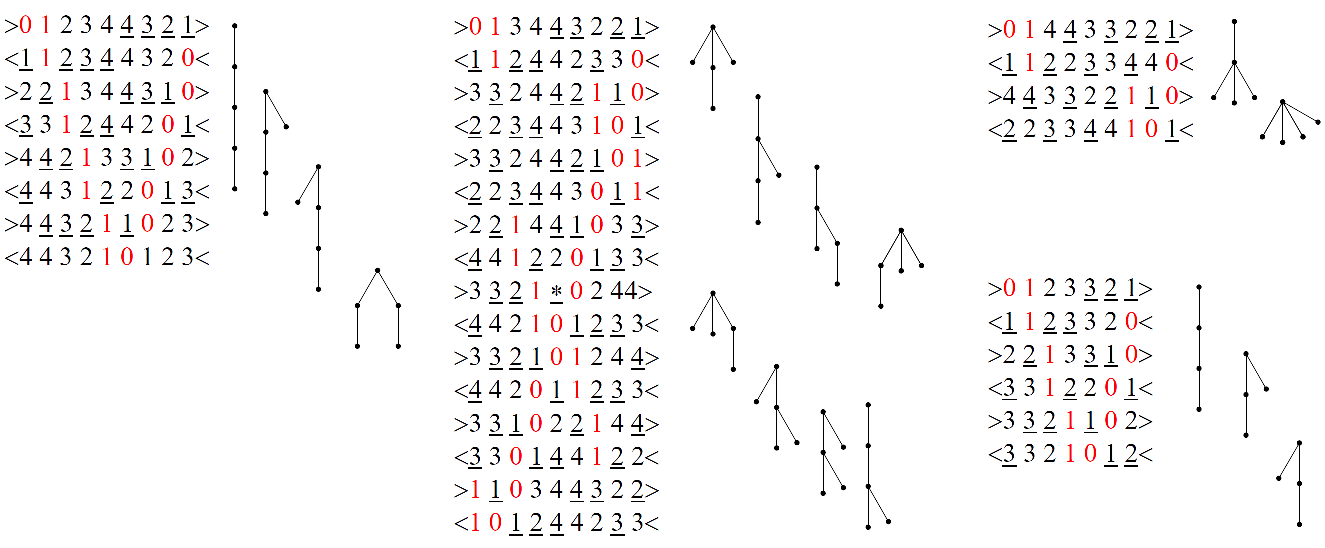}
\caption{Examples of $\lambda=4,8,2,3$.}
\label{gris}
\end{figure}

\begin{remark}\label{remark} Our definition of $f$ (differing from the $f$ in \cite[display (6)]{mmm}) is illustrated in Figure~\ref{gris}, below, in the notation of \cite{D1}, arising from the Dyck path associated to each vertex $x$ of $M_n$, with 0, or the first (resp. second) appearance of each integer in $[n+1]$, corresponding to a 0- (resp. 1-) bit in $x$, where vertices $x\in A_n$ (resp. $x\in B_n$) are expressed as ``$>\ldots >$" (resp. ``$<\ldots <$''), to be read from left to right (resp. right to left).  The ordered trees for those vertices $x$ are drawn to the right of each resulting periodic path in Figure~\ref{gris}.
\end{remark}

\begin{lemma}\cite[Lemma 6]{mmm} Let $n\ge 1$ and let $x\in V(M_n)$. Then, 
\begin{enumerate}
\item $\forall y\in\langle x\rangle$ and $\forall 0\le i\in\mathbb{Z}$, $\langle f^i(x)\rangle=\langle f^i(y)\rangle$. In particular, $\kappa(y)=\kappa(x)$.
\item $\forall 0\le i\in\mathbb{Z}$, $\langle f^i(x)\rangle=\langle f^{\kappa(x)+i}(x)\rangle$.
\item $\forall 0\le i<j\le\kappa(x)$ in $\mathbb{Z}$, $\langle f^i(x)\rangle\ne\langle f^j(x)\rangle$.
\item $\forall 0\le i\in\mathbb{Z}$, $\kappa(f^i(x))=\kappa(x)$.
\item $\kappa(x)=2\lambda(t(x))$, so $\lambda(t(x))$ is semi-period of $f$ at $x$.
\end{enumerate}\end{lemma}

\section{Periodic paths and gluing pairs}\label{per}  

\noindent{\bf Cycle factor of $N_n$:} For any $y\in\langle x\rangle$ and any $0\le i\in\mathbb{Z}$, we have $\kappa(f^i(y))=\kappa(x)$, so $\langle P(x)\rangle=\langle P(f^i(y))\rangle$. This yields a 2-factor of $N_n$, denoted ${\mathcal F}:=\{\langle P(x);x\in V(M_n)\}$.

\begin{proposition}\label{mathcalf}\cite[Proposition 7]{mmm}
For any $n\ge 2$, ${\mathcal F}_n$ has the following properties:
\begin{enumerate}\item
for every $x\in V(M_n)$ the $(2i)$-th vertex $y$ after $x$ on $P(x)$ satisfies $t(y)=\rho^i(t(x))$. Therefore, both $P(x)$ and $\langle P(x)\rangle$ can be identified with $[t(x)]$.
\item $|V(P(x))|=2\lambda(t(x))\ge 4$ and $\ell(f^{2i}(x))=\ell(x)+i$, $\forall i=0,\ldots,\lambda(t(x))$.
\item The cycles of ${\mathcal F}_n$ are in bijective correspondence with the plane trees on $n+1$ vertices.
\end{enumerate}
\end{proposition}

By Proposition~\ref{mathcalf} item 3, the number of cycles of ${\mathcal F}_n$ fits the sequence \cite[\seqnum{A002995}]{oeis}. Also, \cite{mmm} mentions that the number of plane trees, or cycles of ${\mathcal F}_n$, grows exponentially.\\

\noindent{\bf Gluing pairs:} Consider the {\it star} $s_n=0(01)^{n-1}1\in D_n$ for $n\ge 3$ and the {\it footed-star} $s'_n=01s_{n-1}\in D_n$ for $n\ge 4$. A {\it gluing pair} is a pair $(x,y)\ne(s_n,s'_n)$, with $x=u0v011$ and $y=u0v101$, where $u,v\in D$. \\ 

 \noindent{\bf Pull/push operations:} Let $G_n$ be the set of all gluing pairs $(x,y)$, where $x,y\in D_n$. By seeing $x,y$ as rooted trees, \cite{mmm} declares that $y$ is obtained from $x$ by a {\it pull} operation, and calls its inverse a {\it push} operation. In our treatment, the right half of our Figure~\ref{unalinea} resembles, but differs reflectively, from \cite[Figure 9]{mmm}. We write $y=\mbox{pull}(x)$ and $x=\mbox{push}(y)$, say that $x$ is {\it pullable}, $y$ is {\it pushable} and that $u$ and $v$ are the {\it left} and {\it right} subtrees of both $x$ and $y$. 
 
\begin{lemma}\cite[Lemma 8]{mmm}
Let $(x,y)\in G_n$. If $x$ has a centroid in $u$, then $u$ is also a centroid of $y$, so $\phi(y)=\phi(x)-1$. If $y$ has a centroid in $v$, then $v$ is also a centroid of $x$, so $\phi(y)=\phi(x)+1$. 
\end{lemma} 
 
Let $(x,y)\in G_n$. Let $x^i:=f^i(0x)$ and $y^i:=f^i(0y)$, for $i\ge 0$. The resulting sequences agree with the first vertices of $P(x)$ and $P(y)$, respectively.
In such notation, we notice the 6-cycle $C(x,y)=(x^0,y^1,y^0,x^5,x^6,x^1)$, where:

\begin{eqnarray}\label{6}\begin{array}{c}
x^0=0u0v011,\\
y^1=0u0v111,\\
y^0=0u0v101,\\
x^5=0u1v101,\\
x^6=0u1v001,\\
x^1=0u1v011.\\
\end{array}\end{eqnarray}
Then, $P(x)$ and $P(y)$ are glued together by removing the alternate edges $(y^0,y^1)$, $(x^0,x^1)$ and $(x^5,x^6)$ via the symmetric difference between $C(x,y)$ and $P(x)\cup P(y)$.

\begin{lemma}\label{9}\cite[Lemma 9]{mmm}
$(x,y)\in G_n\Rightarrow(|P(x^0)|=\kappa(x^0)\ge 8\mbox{ and }|P(y^0)|=\kappa(y^0)\ge 4)$. 
\end{lemma}

\begin{figure}[htp]
\includegraphics[scale=0.7]{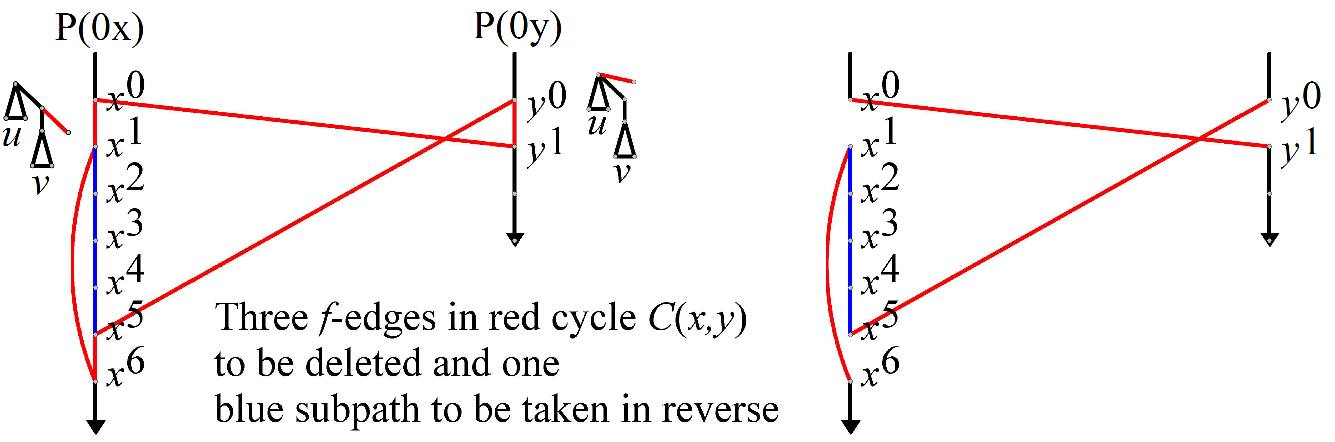}
\caption{The $\nabla$ operation illustrated}
\label{gluing}
\end{figure}

Modifying the proof of Lemma 9 in~\cite{mmm}, we get for our Lemma~\ref{9} now that $$\alpha(C(x,y))=(|u|+|v|+3,|u|+|v|+4,|u|+2,|u|+|v|+3,|u|+|v|+4,|u|+2).$$ 
On the other hand, if $(x,y)=(s_n,s'_n)$, then $\kappa(x^0)=4$, so $\langle x^0\rangle = \langle x^4\rangle$, $\langle x^2\rangle = \langle x^6\rangle$, and Lemma~\ref{9} does not hold.

\begin{remark}\label{fedges} By Lemma~\ref{9}, $\sigma^i(C(x,y))$ shares $\sigma^i(x^0,x^1)$ and $\sigma^i(x^5,x^6)$ with $\sigma^i(P(x^0))$, and $\sigma^i(y^0,y^1)$ with $\sigma^i(P(y^0))$. These edges are said to be the $f${\it-edges} of the {\it gluing cycle} $\sigma^i(C(i,j))$.\end{remark}

If $[x]\ne[y]$, then $\langle P(x)\rangle$ and $\langle P(y)\rangle$ are distinct cycles in $N_n$ by Proposition~\ref{mathcalf}, so we have that the {\it grafted path}
$$P(x^0)\nabla P(y^0):=(x^0,y^1,y^2,\ldots,y^{2\lambda(y)-1},\sigma^{-\lambda(y)}(y^0,x^5,x^4,x^3,x^2,x^1,x^6,x^7,\ldots,x^{2\lambda(x)-1}))$$ is a periodic path in $M_n$.

The $2n+1$ periodic paths 
$\sigma^i(P(x^0)\nabla P(y^0))$ form $$\cup_{i\ge0}(\sigma^i(P(x^0)\nabla P(y^0))),\mbox{ that visits all vertices of }\cup_{i\ge 0}(\sigma^i(P(x^0)\cup P(y^0))).$$ 
Indeed, $|P(x^0)|=2\lambda(x)$, $|P(y^0)|= 2\lambda(y)$ and $\sigma^{\lambda(i)}(y^{2\lambda(y)})=y^0$, by Proposition~\ref{mathcalf}, item 2.
Then, $E(\cup_{i\ge0}(\sigma^i(P(x^0)\nabla P(y^0))))$ is the symmetric difference of $$E(\cup_{i\ge 0}(\sigma^i(P(x^0)\cup P(y^0))))\mbox{ with the gluing cycles }\cup_{i\ge 0}\sigma^i(C(x,y)).$$

\begin{figure}[htp]
\includegraphics[scale=0.7]{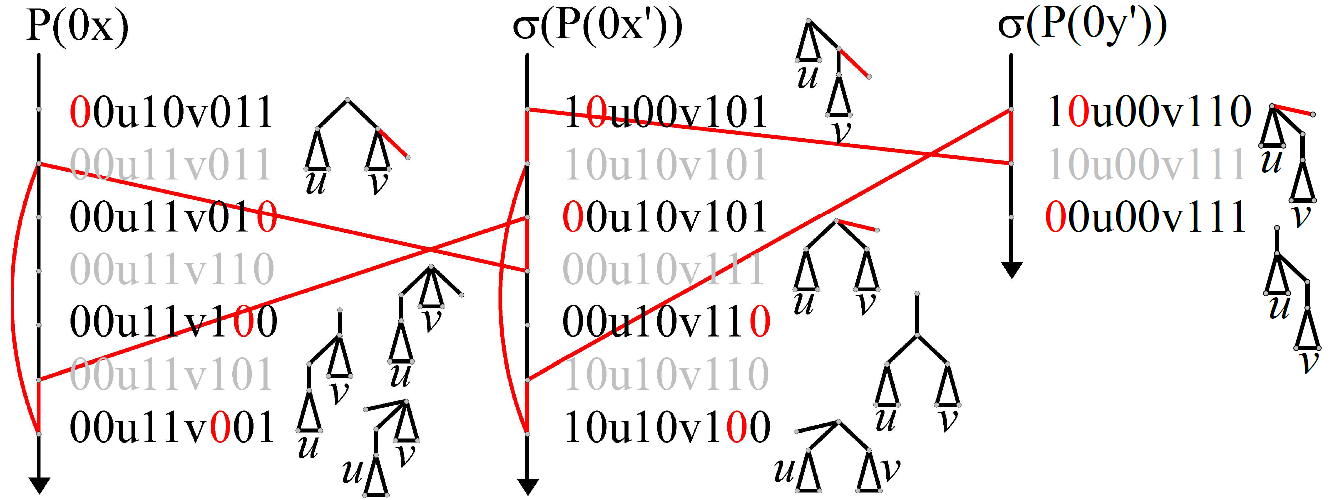}
\caption{Two nested 6-cycles $C(x,y)$ and $C(x',y')$}
\label{achicar}
\end{figure}

\noindent{\bf Additional notation:} For all $i\ge0$, the subpath $\sigma^i(x^1,\ldots,\sigma^5)$ of $\sigma^i(P(x^0))$ is said to be {\it reversed} by $\sigma^i(C(x,y))$. Two gluing cycles $\sigma^i(C(x,y))$ and $\sigma^j(C(x',y'))$ are {\it compatible} if they have no $f$-edges in common. They are {\it nested} if the edge $\sigma^i(y^0,y^1))$ of $\sigma^i(C(x,y))$ belongs to the path reversed by $\sigma^j(C(x',y'))$ (see Figure~\ref{achicar}). They are {\it interleaved} if the $f$-edge $\sigma^j(x'^0,x'^1)$ of $\sigma^j(C(x',y'))$ belongs to the path that is reversed by $\sigma^i(C(i,j))$.  

\begin{proposition}\cite[Proposition 10]{mmm}\label{ojo}
Let $n\ge 4$. Let $(x,y),(x'y')\in G_n$ with $[x]\ne[y]$, $[x']\ne[y']$ and $\{[x],[y]\}\ne\{[x'],[y']\}$. Then, $\forall 0\le i,j\in\mathbb{Z}$, $\sigma^i(C(x,y))$ and $\sigma^j(C(x',y'))$ are:
\begin{enumerate}
\item compatible;
\item interleaving $\Leftrightarrow$ $i=j+2$ and $x'=\rho^2(x)$;
\item nested $\Leftrightarrow$ $i=j-1$ and $x'=\rho^{-1}(y)$.
\end{enumerate}
\end{proposition}

\noindent Item 3 here can be interpreted as follows: Starting at the tree $x$, pull an edge $e$ towards the root to reach the tree $y=$pull$(x)$, then perform an inverse tree rotation $x'=\rho^{-1}(y)$
which makes $e$ pullable, and pull it again to reach $y'=$pull$(x')$. Thus, nested gluing cycles occur if and only if the same edge of the underlying plane trees is pulled twice  in succession.\\

\begin{definition} For $n\ge 4$, let ${\mathcal H}_n$ be the directed arc-labeled multigraph with vertex set $PT_n$ and such that for each $(x,y)\in G_n$ there is an arc labeled $(x,y)$ from $[x]$ to $[y]$. 
\end{definition}

Some pairs of nodes in ${\mathcal H}_n$ may be connected by multiple arcs similarly oriented but with different labels, e.g.
$([0011001101],[0101001101])$ and $([0011010011],[0101010011])$; oppositely oriented, e.g. $([00101011],[01001011])$ and
 $([00110101],[01010101])$. There may be also loops in ${\mathcal H}_n$, e.g. $([00101101],[01001101])$.

\begin{figure}[htp]
\includegraphics[scale=0.7]{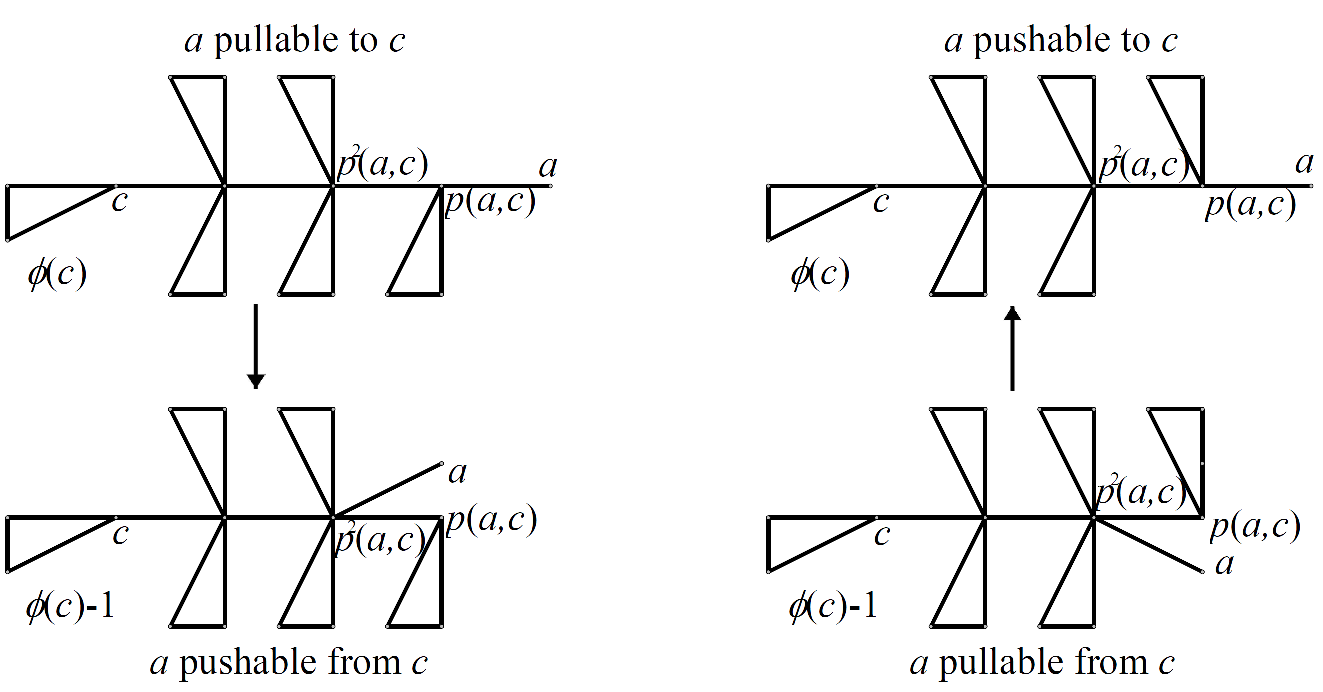}
\caption{pullable and pushable trees}
\label{clokwise}
\end{figure}

\begin{remark}\label{r2} Let $\mathcal T$ be a simple subgraph of ${\mathcal H}_n$. Let $G({\mathcal T})$ be the set of all arc labels of $\mathcal T$. Since $\mathcal T$ is simple, then $[x]\ne[y]$,
$[x']\ne[y']$ and $\{[x],[y]\}\ne\{[x'],[y']\}$, for all $([x],[y]),([x'],[u'])\in G({\mathcal T})$. We say that $G({\mathcal T})$ is {\it interleaving-free} or {\it nesting-free}, respectively, if there are no two gluing pairs $(x,y),(x',y')\in G({\mathcal T})$ such that the gluing cycles $\sigma^i(C(x,y))$ and $\sigma^j(C(x',y'))$ are interleaved or nested for any $i,j\ge 0$.
\end{remark}

\begin{figure}[htp]
\includegraphics[scale=0.87]{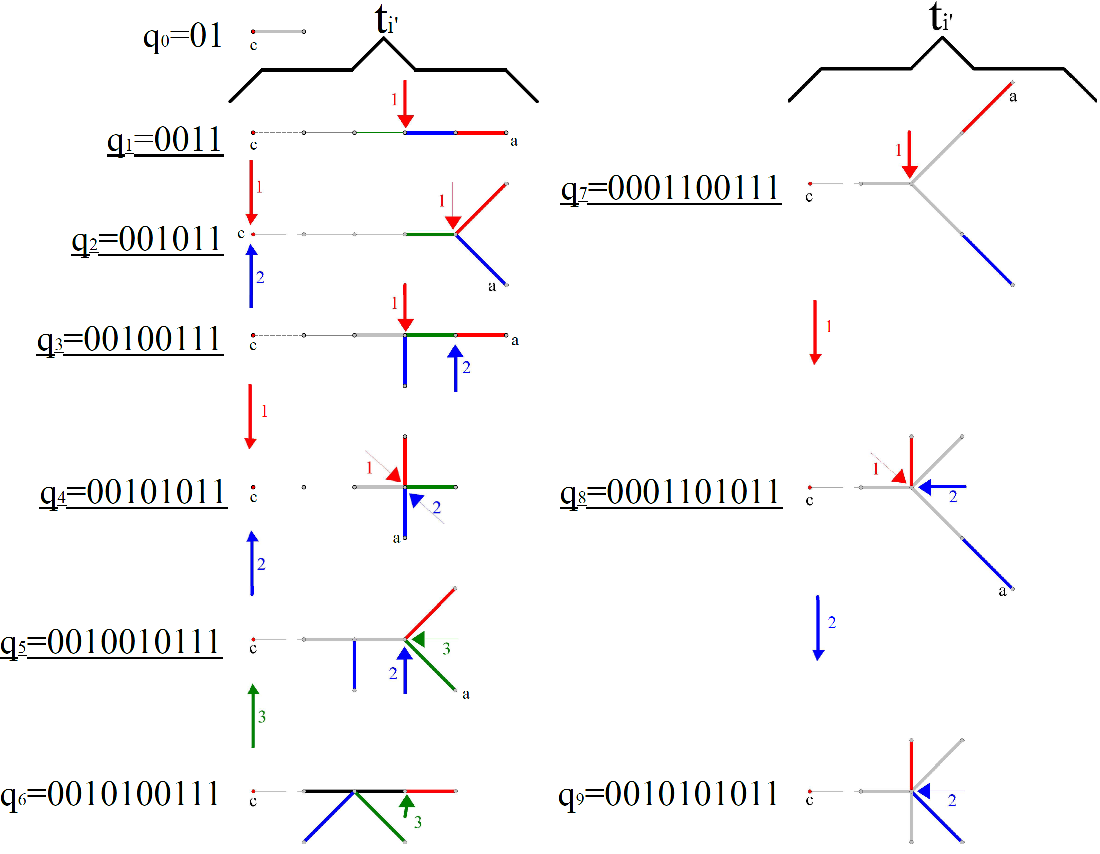}
\caption{Illustration of the trees $q_0,\ldots,q_9$ in thick trace. Remaining edges in path to selected centroid $c$ in thin/dashed trace. Pull/push operations indicated in non-gray/black colors and numbers. The spanning tree ${\mathcal T}_n$ has every arc $([x],[y])$ labeled by a gluing pair $(x,y)$. The rooted trees $x$ and $y$ are obtained by rooting the plane trees $[x]$ and $[y]$ as indicated by the short arrows, thus showing the splitting of the cyclic ordering of the neighbors of each such vertex to get the right-to-left ordering of the children of the corresponding root. Moreover, the short arrow at an $[x]$ has a thick shaft and that of the corresponding $[y]$ a thin shaft. Every arc and the corresponding two short arrows are marked by the same integer. The underlined trees $q_1,\ldots q_5$ and $q_7,q_8$ are treated by separate rules in step (T2).}
\label{q0aq9}
\end{figure}

\begin{lemma}\cite[Lemma 11]{mmm}
If for every $(x,y)\in G({\mathcal T})$ the root of $x$ is not a leaf, then $G({\mathcal T})$ is interleaving-free.
\end{lemma}

\noindent{\bf Pullable/pushable trees:} (See Figure~\ref{clokwise}).
Let $T$ be a tree, let $a,c\in V(T)$ and let $d(a,c)$ be the distance between $a$ and $c$. Let $p^i(a,c)$ be the $i$-th vertex in the path from $a$ to $c$, ($i=0,1,\ldots,d(a,c)$).
In particular, $p^0(a,c)=0$ and $p^{d(a,c)}(a,c)=c$.

Let $c,a\in V(T)$, where $a$ is a leaf of $T$ and $d(a,c)\ge 2$. Then $a$ is {\it pullable to} $c$ if $p^1(a,c)$ has no neighbors between $p^2(a,c)$ and $a$ in the clockwise ordering of neighbors.  (This and the next concepts differ from the couterclockwise stance in \cite{mmm}).
Also, $a$ is {\it pushable to} $c$ if $p^1(a,c)$ has no neighbors between $a$ and $p^2(a,c)$ in the clockwise ordering of neighbors.

Let $d(a,c)\ge 1$. Then, $a$ is {\it pullable from} $c$ if $d(a,c)\ge 2$ and $p^1(a,c)$ has at least one neighbor between $p^2(a,c)$ and $a$ in its clockwise ordering of neighbors or if $c$ is not a leaf and $d(a,c)=1$. Also, $a$ is {\it pushable from} $c$ if $d(a,c)\ge 2$ and $p^1(a,c)$ has at least one neighbor between $a$ and $p^2(a,c)$ in its clockwise ordering of neighbors or if $c$ is not a leaf and $d(a,c)=1$.

For odd $n\ge 5$, consider the {\it dumbbells} $d_n:=(01)^{\frac{n-1}{2}}0(01)^{\frac{n-1}{2}}1$ and $d'_n:=\rho^{2}(d_n)=010(01)^{(n-1)/2}0(01)^{(n-3)/2}$. Each dumbbell has two centroids of degree $\frac{n+1}{2}$, while all the remaining vertices are leaves.
 
If $T$ has just one centroid $c$, every $c$-subtree of $T$ is said to be {\it active}. If $T$ has two centroids $c,c'$, every $c$-subtree except the one containing $c'$ and every $c'$-subtree except the one containing $c$ are also said to be {\it active}. For $n\ge 4$, if $T\ne[s_n]$ and $T\ne[d_n]$ for odd $n$, then $T$ has a centroid with an active subtree that is not a single edge.

\begin{lemma}\label{ay}\cite[Lemma 12]{mmm}
Let $c$ be a centroid of a plane tree $T$, let $a$ be a leaf of $T$ that is pullable to $c$ and that belongs to an active $c$-tree unless $n\ge 5$ is odd with $T=d_n$. Then, the rooted tree $x:=x(T,c,a)=T^{(p^2(a,c),p^1(a,c))}$ is a pullable tree, the rooted tree $y:={\rm pull}(x)$ satisfies $\phi(y)=\phi(x)-1$ and the leaf $a$ is pushable from $c$ in $[y]$.
Moreover, the centroids of $x$ and $y$ are identical and contained in the left subtrees of $x$ and $y$, unless $n\ge 5$ is odd with $x=d_n$, in which case $x$ has two centroids, namely the roots of its left and right subtrees, and the root of the left subtree is the unique centroid of $y$. 
\end{lemma}

A leaf of $T$ is {\it thin} if its unique neighbor in $T$ has degree $\le 2$; otherwise, it is {\it thick}.\\

\begin{lemma}\label{ayayay}\cite[Lemma 13]{mmm}
Let $c$ be a centroid of a plane tree $T$, let $a$ be a thick leaf of $T$ that is pushable to $c$ and that belongs to an active $c$-subtree unless $n\ge5$ is odd with $T=[d'_n]$. Then, the rooted tree $y:=y[T,c,a]:=T^{(p^1(a,c),a)}$ is a pushable tree, the rooted tree $x={\rm push}(y)$ satisfies $\phi(x)=\phi(y)-1$, and the leaf $a$ is pushable from $c$ in $[y]$. Moreover, the centroid(s) of $x,y$ (is) are identical and contained in the right subtrees of $x,y$, unless $n\ge 5$ is odd with $x=d_n$, in which case $x$ has two centroids, namely the roots of its left and right subtrees, and the leaf of its right subtree is the unique centroid of $y$.
\end{lemma}

\begin{definition}\label{tn} For $n\ge 4$, let ${\mathcal T}_n$ be a subgraph of ${\mathcal H}_n$ such that: {\bf(a)} for every $T\in PT_n$ with $T\ne[s_n]$, and $T\ne[d_n]$ if  $n$ is odd, there is a centroid $c$ of $T$ with at least one active $c$-subtree $C$ that is not a single edge; the rightmost leaf of every such $C$ is pullable to $c$; we fix one such leaf $a$; {\bf(b)} If $n$ is odd and $T=[d_n]$, let $c$ be one of its centroids with exactly one $c$-subtree $C$ which is not a leaf, namely the tree $s_{(n+1)/2}$; the rightmost leaf of $C$ is pullable to $c$. In both cases, let $x:=x(T,c,a)$ be the corresponding pullable rooted tree as defined in Lemma~\ref{ay} and define $y:=\mbox{pull}(x)$, yielding the gluing pair $(x,y)\in G_n$. We let ${\mathcal T}_n$ be the spanning subgraph of ${\mathcal H}_n$ given by the union of arcs $([x],[y])$ labeled $(x,y)$ for all gluing pairs $(x,y)$ obtained this way. Ties between two centroids or multiple $c$-subtrees are broken arbitrarily.
For any arc $(T,T')$, $T'$ is an {\it out-neighbor} and $T$ is an {\it in-neighbor}.\end{definition}

\begin{lemma}\label{ahora}\cite[Lemma 14]{mmm}
For any $n\ge 4$, ${\mathcal T}_n$ is a spanning tree of ${\mathcal H}_n$, and for every arc $(T,T')$ in ${\mathcal H}_n$, $\phi(T')=\phi(T)-1$. Every plane tree $T\ne[s_n]$ has exactly one neighbor $T'$ in ${\mathcal T}_n$ with $\phi(T')=\phi(T)-1$ which is an out-neighbor. Furthermore, $G({\mathcal T}_n)$ is interleaving-free.
\end{lemma}

\begin{definition}
Consider a periodic path $P=(x_1,\ldots,x_k)$ in $M_n$. An integer sequence $\alpha=(a_1,\ldots,a_k)$ is a {\it flip sequence} if $a_i$ is the position at which $x_{i+1}$ differs from $x_i$, for each $i\in[k-1]$, and the vertex $x_{k+1}$ obtained from $x_k$ by flipping the bit at position $a_k$ satisfies $\langle x_{k+1}\rangle=\langle x_1\rangle$.
There is a unique integer $\lambda$ mod $2n+1$ given by the relation $x_1=\sigma^\lambda(x_{k+1})$. Let $\lambda(\alpha)=\lambda$ be said to be the {\it shift} of $\alpha$.\end{definition}

\section{Initial attempt at proving Theorem~\ref{the}:}\label{1st} 

\noindent {\bf Scaling trick:}\label{sca} \cite{mmm} presents a {\it scaling trick} consisting in the construction of a flip sequence $\alpha_0$ for one particular shift $s$ coprime to $2n+1$. From such a shift $s$, a transformation $\Upsilon$ yields every shift $s'$ coprime to $2n+1$, where $\Upsilon$ consists in multiplying all entries of $\alpha_0$ by $s^{-1}s'$ mod $2n+1$, with $s^{-1}$ equal to the multiplicative inverse of $s$.

 Define $\label{rev}\mbox{rev}(P):=(x_1,\sigma^{\lambda(\alpha)}(x_k,x_{k-1},\ldots,x_2))$ and 
$\mbox{rev}(\alpha):=(a_k,a_{k-1},\ldots,a_1)-\lambda(\alpha)$, with indices taken mod $2n+1$. Note that rev$(\alpha)$ is a flip sequence for the periodic path rev$(P)$ satisfying $\lambda($rev$(\alpha))=-\lambda(\alpha)$.

Define mov$(P):=(x_2,\ldots,x_k,\sigma^{-\lambda(\alpha)}(x_1))$ and mov$(\alpha):=(a_2,\ldots,a_k,a_1+\lambda(\alpha))$, this being a flip sequence for the periodic path mov$(P)$ satisfying $\lambda($mov$(\alpha))=\lambda(\alpha)$, which means that the shift  is independent of the choice of the starting vertex along the path. Similarly, $\alpha +i$ is a flip sequence for $\sigma^{-i}(P)$ satisfying $\lambda(\alpha +i)=\lambda(\alpha)$, $\forall i\in{\mathbb Z}$. 

\begin{figure}[htp]
\includegraphics[scale=0.72]{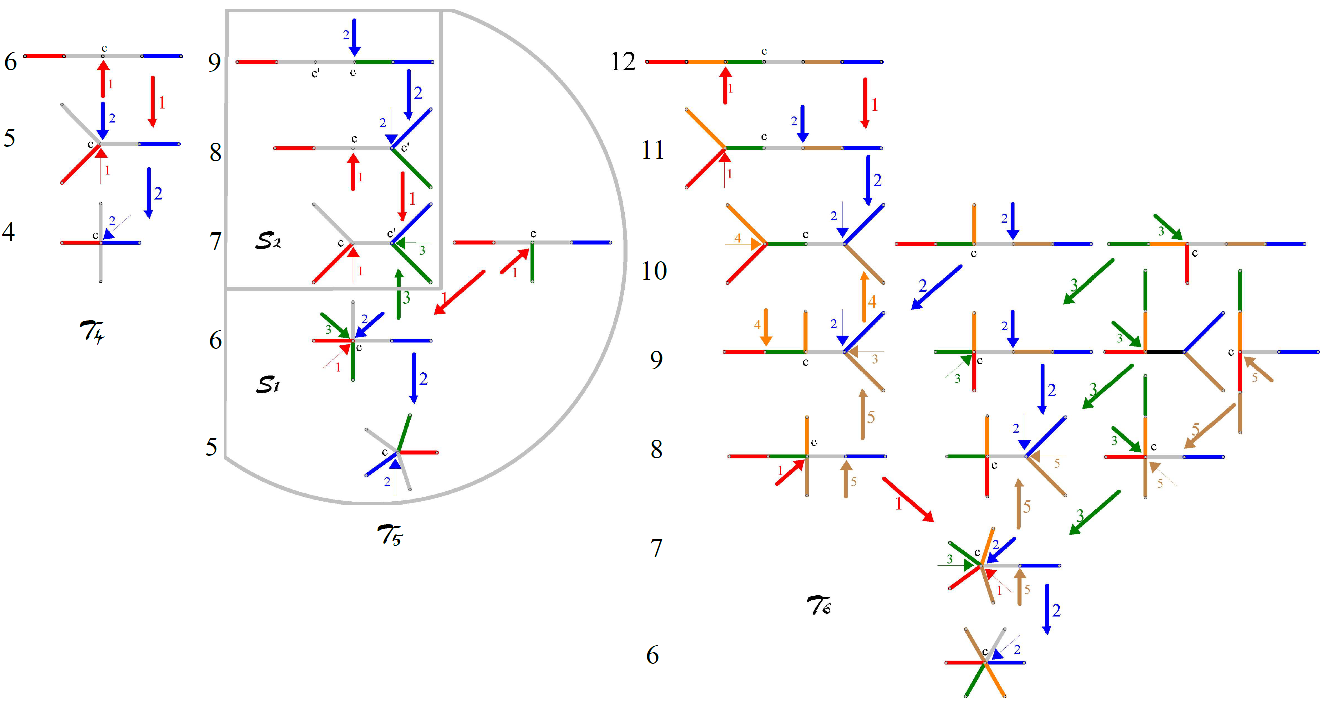}
\caption{Illustration of the spanning trees ${\mathcal T}_4,{\mathcal T}_5,{\mathcal T}_6$. The subgraphs ${\mathcal S}_1,{\mathcal S}_2\subseteq{\mathcal T}_n$ with all plane trees that have one and two centroids, respectively, are highlighted. Centroid(s) are marked red and black, with each red centroid as selected in step (T1). Plane trees are arranged in levels according to their potential, which is shown on the left side. The arrow markings are as in Figure~\ref{q0aq9}.}
\label{t4ts}
\end{figure}

For any $x\in V(M_n)$, let $\alpha(x)$ be the sequence of positions at which $f^{i+1}(x)$ differs from $f^i(x)$, $\forall i=0,\ldots,\kappa(x)-1$. 
Clearly, $\alpha(x)$ is a flip sequence for $P(x)$ (Definition~\ref{p}).
By Proposition~\ref{mathcalf} item 2,
$\lambda(\alpha(x))=\lambda(t(x))$. 

For any subtree ${\mathcal T}$ of ${\mathcal H}_n$ with $G:=G({\mathcal T})$ interleaving-free as in Remark~\ref{r2}, define the set of necklaces 
$N({\mathcal T}):=\cup_{[x]\in{\mathcal T}}\langle P(0x)\rangle$. By Proposition~\ref{mathcalf} item 1, this is the set of necklaces visited by those cycles $\langle P(0x)\rangle$ in $N_n$ for which $[x]\in\mathcal T$.

For any $z\in N({\mathcal T})$ and any $x\in z$, there is a pair ${\mathcal P}_G(X)=\{P,P'\}$ of two periodic paths $P$ and $P'$, both starting at $x\in V(M_n)$ with flip sequences $\alpha(P)$ and $\alpha(P')$ such that $P'=$rev$(P)$ and $\alpha(P')=$rev$(\alpha(P))$. Moreover, $\langle P\rangle$ and $\langle P'\rangle$ are oppositely oriented in the subgraph of $N_n$ whose vertex set is $N({\mathcal T})$.

The node set of ${\mathcal T}_n$ is $PT_n$. By Lemma~\ref{ahora}, $G({\mathcal T}_n)$ is interleaving-free. Fix $x_1:=0^{n+1}1^n\in V(M_n)$. The pair ${\mathcal P}_{G({\mathcal T}_n)}(x_1)$ contains a periodic path $P$ with starting vertex $x_1$ and second vertex $f(x_1)$ in $M_n$ such that $\langle P\rangle$ has vertex set $N({\mathcal T}_n)=\cup_{[x]\in PT_n}\langle P(0x)\rangle=\{\langle x\rangle |x\in V(M_n)\}$, i.e. $\langle P\rangle$ is Hamilton cycle in $N_n$. The corresponding flip sequence $\alpha(P)$ has a shift
\begin{eqnarray}\label{la}\lambda(\alpha(P))=\sum_{T\in PT_n}\gamma_T.\lambda(T),\end{eqnarray} for numbers $\gamma_T\in\{1,-1\}$ that are determined by which gluing cycles encoded by ${\mathcal T}_n$ are nested.

With $s:=\lambda(\alpha(P))$, define $\alpha_0:=\alpha(P)$ and $\alpha_i:=\alpha_0+i.s$, for $i\in[2n]$. If we apply the entire flip sequence $(\alpha_0,\alpha_1,\ldots,\alpha_{2n})$ to $x_1$ in $M_n$, then the vertex $\sigma^{i.s}(x_1)$ is reached after applying all flips in $(\alpha_0,\alpha_1,\ldots,\alpha_{i-1})$, $\forall i\in[2n+1]$. If $s$ and $2n+1$ are coprime, then $x_1$ is reached only after applying the entire flip sequence. Since $\alpha_0$ is the flip sequence of the Hamilton cycle $\langle P\rangle$ in $N_n$, the resulting  sequence of bitstrings is a Hamilton cycle in $M_n$. However, this approach requires that $s=\lambda(\alpha(P))$ and $2n+1$ be coprime. A technique is needed to modify $\alpha(P)$ into another flip sequence $\alpha'$ such that $s':=\lambda(\alpha')$ is coprime to $2n+1$.\\

\noindent{\bf Switches and their shifts:} Let $p(x,y)$ indicate a pair $x,y\in V(M_n)$ that differ in just one position. A triple of vertices $\tau=(x,y,y')$, where $x\in A_n$, $\{y,y'\}\subseteq B_n$ and $y\ne y'$, 
 is a {\it switch} if $x$ differs from $y$, (resp. $y'$) in a single bit, and $\langle y\rangle=\langle y'\rangle$. In $N_n$, a switch may be considered as a multiedge $(\langle x\rangle,\langle y\rangle)=(\langle x\rangle,\langle y'\rangle)$. The {\it shift} of a switch $\tau=(x,y,y')$, denoted $\lambda(\tau)$, is the integer $i$ such that $y=\sigma^i(y')$. Denote a switch $\tau=(x,y,y')$ compactly by writing $x$ with the 0-bit at position $p(x,y)$ underlined and the 0-bit at position $p(x,y')$ overlined. This way, we write $\tau=(0000111,1000111,0001111)=\underline{0}00\overline{0}111$. For any switch $\tau=(x,y,y')$, the inverted switch $\tau^{-1}=(x,y',y)$ has shift $\lambda(\tau^{-1})=-\lambda(\tau)$. Clearly, cyclically rotating a switch  yields another switch with the same shift. Also, reversing a switch yields another switch with the negated shift. For example,
 $\sigma(\tau)=11\overline{0}00\underline{0}1$ has shift 1 while its reversed switch $1\underline{0}00\overline{0}11$ has shift $-1$. 

Consider a flip sequence $\alpha=(a_1,\ldots,a_k)$ with shift $\lambda(\alpha)$ for a periodic path $P=(x_1,\ldots,x_k)$ and let $x_{k+1}$ be the vertex obtained from $x_k$ by flipping the bit at position $a_k$. If $(x_i,x_{i+1})=(x,y)$ for some $i\in[k]$, then the modified flip sequence 
\begin{eqnarray}\label{25a} \alpha'=(a_1,\ldots,a_{i-1},p(x,y'),a_{i+1}+\lambda(\tau),\ldots,a_k+\lambda(\tau))\end{eqnarray}
yields a periodic path $P'=(x'_1,\ldots,x'_k)$ that visits necklaces in the same order as does $P$, i.e. $\langle x_i\rangle=\langle x'_i\rangle$, for $i\in[a_k]$, and $\lambda(\alpha')=\lambda(\alpha)+\lambda(\tau)$.
The situation where $(x_i,x_{i+1})=(x,y')$ is symmetric and considers the inverted switch $\tau^{-1}$, with $\lambda(\tau^{-1})=-\lambda(\tau)$.

Similarly, if $(x_i,x_{i+1})=(y',x)$ for some $i\in[k]$, then the modified sequence 
\begin{eqnarray}\label{25c}\alpha':=(a_1,\ldots,a_{i-1},p(x,y)+\lambda(\tau),a_{i+1}+\lambda(\tau),\ldots,a_k+\lambda(\tau))\end{eqnarray} produces a periodic path $P'=(x'_1,\ldots,x'_k)$ that visits necklaces in the same order as $P$, and $\lambda(\alpha')=\lambda(\alpha)+\lambda(\tau)$. The situation $(x_i,x_{i+1})=(y,x)$ is symmetric and goes to the inverted switch $\tau^{-1}$, with $\lambda(\tau^{-1})=-\lambda(\tau)$.

 In particular, if $\langle P\rangle$ is a Hamilton cycle in $N_n$, then $\langle P'\rangle$ is also a Hamilton cycle in $N_n$ with shift $\lambda(\alpha')=\lambda(\alpha)+\lambda(\tau)$.\\
 
 \noindent{\bf Construction of switches out of $\tau_{n,1}=\underline{0}0^{n-1}\overline{0}1^n$, with $\lambda(\tau)=1$:}
 
 For any integers $n\ge 1$, $d\ge 1$ and $1\le s\le d$, make the $(s,d)${\it -orbit} to be the maximal prefix of the sequence $s+id$, $i\ge0$, modulo $2n+1$, in which all the numbers are distinct. Then, the number of distinct $(s,d)$-orbits for fixed $d$ and $s\ge 1$ is $n_d:=$gcd$(2n+1,d)$, and the length of each orbit is $\ell_d:=\frac{2n+1}{n_d}$, where both $n_d$ and $\ell_d$ are odd. For example, let $n=10$ and $d=6$, so $n_d=3$, $\ell_d=7$ and the $(1,6)$-orbit is $(1,7,13,19,4,10,16)$, the $(2,6)$-orbit is $(2,8,14,20,5,11,17)$ and the $(3,6)$-orbit is $(3,9,15,21,6,12,18)$.
 
 For any integer $d$ ($2\le d\le n$) that is coprime to $2n+1$, let $\tau_{n,d}$ denote the sequence whose entries at the positions given by the $(1,d)$-orbit equal $\tau_{n,1}=\underline{0}0^{n-1}\overline{0}1^n$, including the overlined and underlined entries.
 
For any $n\ge 1$, let $Z_n$ be the set of bitstrings of length $2n$ and weight $n$.
For any integers $d$ ($3\le d\le n$) not coprime to $2n+1$, select an arbitrary bitstring $z=(z_2,\ldots,z_{n_d})\in Z_{(n_d-1)/2}$. Let $\tau_{n,d,z}$ be the sequence: {\bf(a)} whose entries at the positions given by the $(1,d)$-orbit form the sequence $\tau_{(\ell_d-1)/2,1}$, including the underlined and overlined entries, and {\bf(b)} for $j=2,\ldots,n_d$, all entries at the positions given by the $(j,d)$-orbit form the sequence $z_j$. Then, the number of choices for $z$ in such a construction is ${n_d-1\choose (n_d-1)/2}$.

\begin{lemma}\label{ta}\cite[Lemma 15]{mmm}
Let $n\ge 1$. For any integer $d$ ($1\le d\le n$) coprime to $2n+1$, the sequence $\tau_{n,d}$ is a switch with $\lambda(\tau_{n,d})=d$. For any integer $3\le d\le n$ not coprime to $2n+1$ and any bitstring $z\in Z_{(n_d-1)/2}$, the sequence $\tau_{d,n,z}$ is a switch with $\lambda(\tau_{d,n,z})=d$.
\end{lemma}
 
\noindent{\bf Interactions:} Given a flip bijection $f$ of $V(M_n)$ as in Section~\ref{flip}, a switch $\tau=\tau(x,y,y')$ is said to be $f${\it-conformal} if either $y=f(x)$ or $x=f(y')$; in such cases, $(x,y)$ or $(y',x)$, respectively, is said to be the $f$-edge of $\tau$. We say that $\tau$ is $f^{-1}$-{\it conformal} if $\tau^{-1}$ is conformal and we refer to the $f$-edge of $\tau^{-1}$ also as an $f$-edge of $\tau$. That a switch is $f$-conformal means that its $f$-edge belongs to a periodic path, as in Definition~\ref{p}.

Given a subset $G\subseteq G_n$, an $f$-conformal or $f^{-1}$-conformal switch $\tau$ is {\it usable with respect to} $G$ if for every $(x',y')\in G$ and all $i\ge 0$, the three $f$-edges of $\sigma^i(C(x',y'))$ in Remark~\ref{fedges} are distinct from the $f$-edges of $\tau$. Those three $f$-edges are removed when joining periodic paths.

\begin{lemma}\cite[Lemma 16]{mmm}
Let $\tau=(x,y,y')$ be an $f^{-1}$-conformal switch with $f$-edge $(y,x)$ for which $t(x)=00\ldots$. Then $\tau$ is usable  with respect to any subset $G\subseteq G_n$.
\end{lemma}

\begin{lemma}\label{tai}\cite[Lemma 17]{mmm}
Let $n\ge 4$. The switch $\tau_{n,1}=:(x,y,y')=\underline{0}0^{n-1}\overline{0}1^n$, where $x\in A_n $ and $y\in B_n$ differ in the first bit, has $f$-edge $(y,x)$ and is $f^{-1}$-conformal. The switch $\tau_{n,2}=:\overline{0}\underline{0}(01)^{n-1}=(x,y,y')=\overline{0}\underline{0}1(01)^{n-1}$, where $x\in A_n$ and $y'\in B_n$ differ in the first bit, has $f$-edge $(y',x)$ and is $f$-conformal. Both switches are usable with respect to any subset $G\subseteq G_n$.
\end{lemma}
 
\begin{lemma}\label{11}\cite[Lemma 18]{mmm}
Let $n\ge 11$ and let $3\le c,d\in\mathbb{Z}$ be such that $c.d=2n+1$. The switch $\tau_{n,d,z}=:(x,y,y')=z\underline{0}(z0)^{(c-3)/2}z\overline{0}(z1)^{(c-1)/2}$, where $z:=0^{(d-1)/2}1^{(d-1)/2}\in Z_{(d-1)/2}$ has $f$-edge $(y,x)$ and is $f^{-1}$-conformal and usable with respect to  $G({\mathcal T}_n)$.
\end{lemma}

\noindent{\bf Number theory:} Let $n\ge 1$. Let ${\mathcal P}(n)$ be the set of prime factors of $n$. For any $s\in\{0,1,\ldots,n-1\}$, define ${\mathcal P}(n,s):={\mathcal P}(n)\setminus{\mathcal P}(s)$, if $s>0$, and ${\mathcal P}(n,0):=\emptyset$.

\begin{lemma}\label{19}\cite[Lemma 19]{mmm}
Let $n\ge 1$ be such that $2n+1$ is not a prime power. Let $s\in\{0,\ldots,n-1\}$ be not coprime to $2n+1$. If ${\mathcal P}(2n+1,s)\ne\emptyset$, then both numbers $s+d$ and $s-d$, where $d:=\Pi\{p\in{\mathcal P}(2n+1,s)\}$, are coprime to $2n+1$. If ${\mathcal P}(2n+1,s)=\emptyset$, then ${\mathcal P}(2n+1,s+d)={\mathcal P}(2n+1,s-d)={\mathcal P}(2n+1))\setminus\{d\}\ne\emptyset$, for any $d\in{\mathcal P}(2n+1)$.
\end{lemma}

\begin{lemma}\label{both}\cite[Lemma 20]{mmm}
Let $n$ be an integer such that $4\le n\le 10$ and let $s\in\{0,\ldots,2n\}$ be not coprime to $2n+1$. Then, both numbers in at least one of the pairs $\{s-1,s+1\}$, $(s-2,s+2\}$, $\{s-1,s+2\}$, $\{s+1,s-2\}$ are coprime to $2n+1$.
\end{lemma}

\section{Sketch of proof of Theorem~\ref{the}}\label{sketch1} 

\begin{proof} Theorem~\ref{the} is established via the {\it scaling trick} (Section~\ref{1st}) for each $n\ge 1$ and any value of $s$ coprime to $2n+1$. For $n=1$, via flip sequence $\alpha:=32$ starting at $x_1=001$ and yielding shift $s=-1$; for $n=2$, via flip sequence $\alpha=1531$ starting at $00011$ and yielding shift $s=1$; for $k=3$, via flip sequence $2635426753$ starting at $0000111$ and yielding $s=-1$.

Assume $n\ge 4$. Consider the spanning tree $ {\mathcal T}_n\subseteq{\mathcal H}_n$ (Definition~\ref{tn}). In Section \ref{1st}, a periodic path $P$ is defined with starting vertex $x_1=0^{n+1}1^n$ and second vertex $f(x_1)$ in $M_n$ such that $\langle P\rangle$ is a Hamilton cycle in $N_n$ and the shift of the corresponding flip sequence $\alpha(P)$ is given by (\ref{la}). Denote this shift by $s:=\lambda(\alpha(P))$. If $s$ is coprime to $2n+1$, we are done. Let us consider the case $s$ not coprime to $2n+1$.

If $4\le n\le 10$, consider the switches $\tau_{n,1}$ and $\tau_{n,2}$, which are $f^{-1}$- and $f$-conformal, respectively, both usable with respect to $G({\mathcal T}_n)$ by Lemma~\ref{tai}. By Lemma~\ref{ta}, their shifts are $\lambda(\tau_{n,1})=1$ and $\lambda(\tau_{n,2})=2$,  respectively. Consequently, by modifying the flip sequence $\alpha(P)$ to become like $\alpha'$ in (\ref{25a}) via one of the two, or both, switches, we obtain a flip sequence $\alpha'$ with shift \begin{eqnarray}\label{28a}s':=\lambda(\alpha')=s + \chi_1.\gamma_1 + \chi_2.\gamma_2.2,\end{eqnarray} for certain signs $\gamma_1,\gamma_2\in\{1,-1\}$ and with indicators  $\chi_1,\chi_2\in\{0,1\}$ that are nonzero if the corresponding switches are employed.

If $n\ge 11$, we distinguish 3 cases: If $2n+1$ is prime power, then $s$ is also power of the same prime. Apply the switch $\tau_{n,1}$ as above, modifying $\alpha(P)$ so that the resulting flip sequence $\alpha'$ has shift \begin{eqnarray}\label{28b}s':=\lambda(\alpha')+\gamma_1.1,\end{eqnarray} for some $\gamma_1\in\{-1,1\}$, with $s'=s\pm 1$ coprime to $2n+1$.

If $2n+1$ is not a prime power and ${\mathcal P}(2n+1,s)\neq\emptyset$, we define $d:=\Pi\{p\in{\mathcal P}(2n+1,s)\}$ and $c:=\frac{2n+1}{d}$, and consider the switch $\tau_{n,d,z}$ from Lemma~\ref{11}, modifying the flip sequence $\alpha(P)$ so that the resulting $\alpha'$  has shift \begin{eqnarray}\label{28c}s':=\lambda(\alpha')=s+\gamma_d.d.\end{eqnarray}

If $2n+1$ is not a prime power and ${\mathcal P}(2n+1,s)=\emptyset$, then we pick some $d\in{\mathcal P}(2n+1)$, define $c=\frac{2n+1}{d}$ and apply switch $\tau_{n,d,z}$, yielding a flip sequence $\alpha'$ with shift $s'$ given as above, which satisfies ${\mathcal P}(2n+1,z)\neq\emptyset$ by the last sentence of Lemma~\ref{19}. We then modify the flip sequence a second time  as in the previous case, and the switch used is distinct  from the first one, as $d':=\Pi\{p\in{\mathcal P}(2n+1,s\pm d)\}=\Pi\{p\in{\mathcal P} (2n+1)\setminus\{d\}\}$ clearly satisfies $d'\neq d$.\end{proof}

The switches $\tau_{n,1}$ and $\tau_{n,2}$ are not sufficient alone for the proof of Theorem~\ref{the}, starting with $n=52$ and $2n+1=105=3\dot 5\dot 7$ and $s=5$, for which none of the three numbers $s-2=3$, $s+1=6$ and $s+2=7$ is coprime to $2n+1$, so Lemma~\ref{both} cannot be applied and it is necessary  to use switches $\tau_{n,s,z}$.

\begin{remark} If we knew that $G({\mathcal T}_n)$ is not only interleaving-free, but also nesting-free, then this would guarantee that all signs $\gamma_T$ in (\ref{la}) are positive, yielding
\begin{eqnarray}\label{Cat} s=\lambda(\alpha(P))=\sum_{T\in{\mathcal T}_m}\lambda(T)=C_n.\end{eqnarray}
Following \cite{mmm} with the stance of \cite{D1}, we define now another spanning tree ${\mathcal T}_n$ of ${\mathcal H}_n$ such that $G({\mathcal T}_n)$ is both interleaving-free and nesting-free.\end{remark}

\section{Efficient computation: Redefinition of ${\mathcal T}_n$}\label{redef} 

In \cite{mmm}, ten rooted trees are distinguished, that in our alternate viewpoint are expressed as:
\begin{eqnarray}\label{eqn}\begin{array}{llll}
q_0:=01,&q_1:=0011,&q_2:=001011,&q_3:=00100111,\\
q_4:=00101011,&q_5:=0010010111,&q_6:=0010100111,&q_7:=0001100111,\\
q_8:=00011010111,&q_9:=0010101011.&&
\end{array}\end{eqnarray}

 For $n\ge 4$, let ${\mathcal T}_n$ be a subgraph of ${\mathcal H}_n$ given as follows: For each plane tree $T\in PT_n$ with $T\ne[s_n]$, consider a gluing pair $(x,y)\in G_n$ with either $T=[x]$ or $T=[y]$. 
  Let ${\mathcal T}_n$ be the spanning subgraph of ${\mathcal H}_n$ given by the union of the arcs $([x],[y])$, labeled $(x,y)$, for all gluing pairs $(x,y)$ obtained this way. 
 The definition of the gluing pair $(x,y)$ for a given plane tree $T\ne[s_n]$ proceeds as in the following three steps (T1)-(T3), unless $n$ is odd and $T=[d_n]$, in which case the special rule (D) is applied.\\

{\bf (D) Dumbbell rule:} If $n$ is odd and $T=[d_n]$, let $c$ be one of the centroids of $T$ having exactly one $c$-subtree that is not a single edge, namely the subtree $s_{(n+1)/2}$. Its leftmost leaf $a$ is thick and pushable to $c$ in $T$, so we define $y:=y(T,c,a)=d'_n$ and $x:=push(y)$, as in Lemma~\ref{ayayay}. \\

{\bf (T1) Fix the centroid and subtree ordering:} If $T$ has two centroids, we let $c$ be a centroid of $T$ whose active $c$-subtrees are not all single edges. If this is true for both centroids, we let $c$ be the one for which its clockwise-ordered active $c$-subtrees $t_1,\ldots,t_k$ give the lexicographically minimal string $(t_1,\ldots,t_k)$, with $t_1$ as the first $c$-subtree found after the $c$-subtree containing the other centroid.

If $T$ has a unique centroid, we denote it by $c$ and consider all $c$-subtrees of $T$, denoting them $t_1,\ldots,t_k$, i.e. $T=[(t_1,\ldots,t_k)]$ such that among all possible  clockwise orderings around $c$, the string $(t_1,\ldots,t_k)$ is lexicographically minimal.\\

{\bf (T2) Select $c$-subtree of $T$:} If $T$ has two centroids, we let $t_{i'}$ be the first of the trees $t_1,\ldots,t_k$ that differs from $q_0$, in display (\ref{eqn}). 

If $T$ has a unique centroid, then for each of the following conditions (i)-(iv), 
we consider all trees $t_i$ ($i\in[k]$) and determine the first subtree $t_i$ satisfying the condition, 
i.e. only check one of these conditions once all trees failed all previous conditions:

\begin{enumerate}
\item[(i)] $t_i=q_1$ and $t_{i-1}=q_0$;
\item[(ii)] $t_i\in\{q_2,q_4\}$ and $t_{i+1}\in\{q_0,q_1,q_2\}$;
\item[(iii)] $t_i\notin\{q_0,q_1,q_2,q_4\}$;
\item[(iv)] $t_i\neq q_0$.
\end{enumerate}
Conditions (i) and (ii) refer to the previous tree $t_{i-1}$ and the next tree $t_{i+1}$ in the clockwise ordering of $c$-subtrees, and those indices are considered modulo $k$. Note that $T\ne[s_n]$. Thus, at least one $c$-subtree is not $q_0$ and satisfies condition (iv), so this rule to determine $t_i$  is well defined. Let $t_{i'}$ be the $c$-subtree determined in this way. Clearly, $t_{i'}$ has at least two edges.\\

{\bf (T3) Select a leaf to push/pull:} If $t_{i'}=0^lq_j1^l$, for some $l\ge 0$ and $j\in\{1,\ldots,8\}$, i.e. $t_{i'}$ is a path with one of the trees $q_j$ attached to it. Then, four cases are distinguished:
\begin{enumerate}
\item[(q137)] If $j\in\{1,3,7\}$, let $a$ be the rightmost leaf of $t_{i'}$, which is thin, and define $x:=x(T,c,a)$ and 
$y:=\mbox{pull}(x)$, as in Lemma~\ref{ay}. 
Clearly, $(j=1)\Rightarrow(y=0^{l-1}q_21^{l-1}$, if $l>0$, and $y=q_0^2$, if $l=0)$; also, $(j=3)\Rightarrow (y=0^lq_41^l)$; and $(j=7)\Rightarrow (t=0^lq_81^l)$.
\item[(q24)] If $j\in\{2,4\}$, let $a$ be the leftmost leaf of $t_{i'}$, which is thick, and define $y:=y(T,c,a)$ and $x:=\mbox{push}(y)$, as in Lemma~\ref{ayayay}. Clearly, $(j=2)\Rightarrow(x=0^{l-1}q_31^{l-1}$, if $l>0$, and $x=q_2q_0$, if $l=0)$; also, $(j=4)\Rightarrow(x=0^{l-1}q_51^{l-1}$, if $l>0$, and $x=q_2q_0$, if $l=0)$. 
\item[(q5)] If $j=5$, let $a$ be the unique leaf of $t_{i'}$ that is neither the rightmost nor the leftmost leaf of $t_{i'}$, where $a$ is thick, and define $y:=y(T,c,a)$ and $y=\mbox{push}(x)$, as in Lemma~\ref{ayayay}. Clearly, $x=0^lq_61^l$.
\item[(q8)] If $j=8$, let $a$ be the leftmost leaf of $t_{i'}$, which is thin, and define $x:=x(T,c,a)$ and $y:=\mbox{pull}(x)$, as in Lemma~\ref{ay}. Clearly, $y=0^lq_01^l$. 
\end{enumerate}
Otherwise, there are three cases to be distinguished:
\begin{enumerate}
\item[(e)] If the potential $\phi(T)=\phi(c)$ is even, we let $a$ be the rightmost leaf of $t_{i'}$ and define $x:=x(T,c,a)$ and $y=\mbox{pull}(x)$, as in Lemma~\ref{ay}.
\item[(o1)] If the potential $\phi(T)=\phi(c)$ is odd and the leftmost leaf $a$ of $t_{i'}$ is thin, define $x:=x(T,c,a)$ and $y:=\mbox{pull}(y)$, as in Lemma~\ref{ay}. 
\item[(o2)] If the potential $\phi(T)=\phi(c)$ is odd and the leftmost leaf $a$ of $t_{i'}$ is thick, define $y:=y(T,c,a)$ and $x:=\mbox{push}(y)$, as in Lemma~\ref{ayayay}.
\end{enumerate} 
This completes the definition of ${\mathcal T}_n$. 
\cite{mmm} refers to rules (q137), (q8), (e) and (o1) in step (T3) as {\it pull rules} and to rules (q24), (q5) and (o2) as {\it push rules}. The leaf to which one of the pull rules (q137), (q8) or (o1) is applied is always thin, whereas the leaf to which any push rule is applied is always thick.

\begin{figure}[htp]
\includegraphics[scale=1.15]{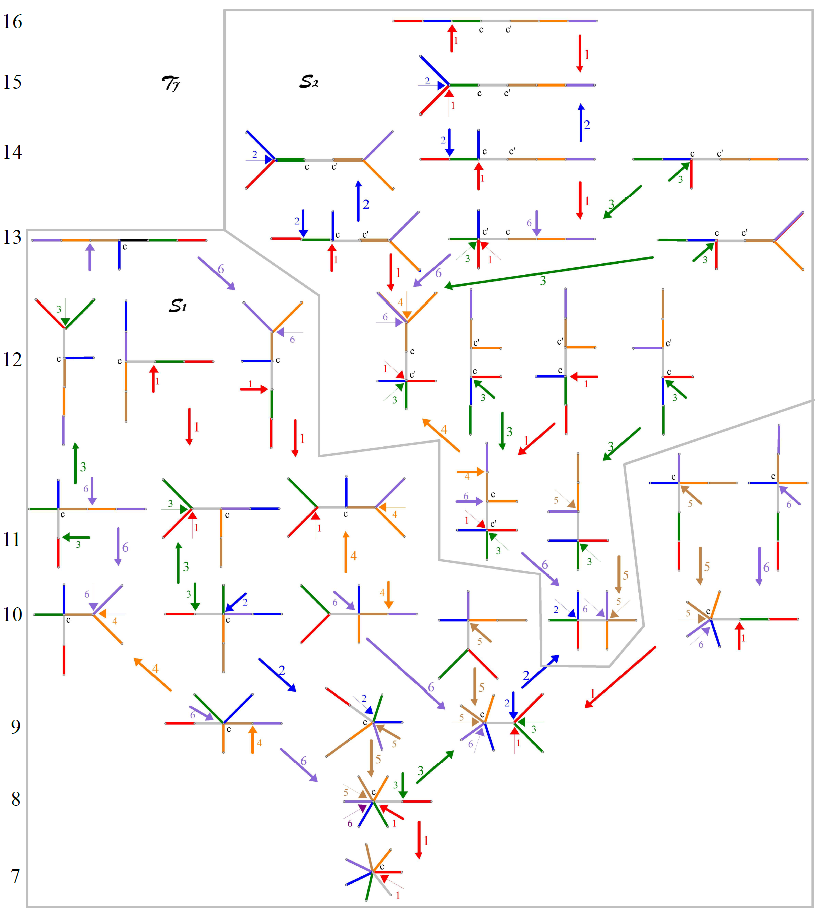}
\caption{Illustration of the spanning tree ${\mathcal T}_7$.}
\label{t4ts'}
\end{figure}

\begin{lemma}\cite[Lemma 21]{mmm} If $T$ has a unique centroid $c$, then the $c$-subtree $t_{i'}$ selected in Step (T2) satisfies the following conditions:
\begin{enumerate}
\item If $t_{i'}=q_1$, then $t_{i'-1}=q_0$ or $t_1=t_2=\cdots=t_k=q_1$.
\item If $t_{i'}=\in\{q_2,q_4\}$, then $t_{i'+1}\in\{q_0,q_1,q_2\}$ or $t_1=t_2=\cdots=t_k=q_4$.
\end{enumerate} 
\end{lemma}

\begin{lemma}\cite[Lemma 22]{mmm}
For any $n\ge 4$, the graph ${\mathcal T}_n$ is a spanning tree of ${\mathcal H}_n$. For every arc $(T,T')$ in ${\mathcal T}_n$, either $\phi(T')=\phi(T)- 1$ or $\phi(T)=\phi(T')- 1$. Every plane tree $T\ne[s_n]$ has exactly one neighbor $T'$ in ${\mathcal T}_n$ with $\phi(T')=\phi(T)-1$ which is either an out-neighbor or in-neighbor. Moreover, $G({\mathcal T}_n)$ is interleaving-free and nesting-free. 
\end{lemma}

Illustrations for the spanning trees ${\mathcal T}_i$ for $i=4,5,6,7$ are in Figures~\ref{t4ts}-\ref{t4ts'}, the versions of~\cite[Figures 15-16]{mmm} for the alternate viewpoint in this survey. \\

\noindent{\bf Interaction with switches:} Assume $G\subseteq G_n$ is nesting-free. A usable switch $\tau$ is {\it reversed} if the $f$-edge of $\tau$ lies on the reversed path of one of the gluing cycles $\sigma^i(C(x',y'))$, $(x',y')\in G$, for some $i\ge 0$, i.e. on the path $\sigma^i(x'^1,\ldots,x'^5)$.

\begin{lemma}\label{ult}\cite[Lemma 23]{mmm} Let ${\mathcal T}_n$ be the spanning tree of ${\mathcal H}_n$ of Section~\ref{redef}. For $n\ge 4$, the switch $\tau_{n,1}$ is reversed with respect to $G({\mathcal T}_n)$ and the switch $\tau_{n,2}$ is not reversed with respect to any set of gluing pairs $G\subseteq G_n$. For $n,d,z$ as in Lemma~\ref{11}, the switch $\tau_{n,d,z}$ is usable and not reversed with respect to $G({\mathcal T}_n)$.
\end{lemma}

\section{Sketch of proof of algorithmic Theorem~\ref{extra}}\label{sketch2}

The algorithm in Theorem~\ref{extra} is a faithful implementation of the constructive proof of Theorem~\ref{the} sketched in Section~\ref{sketch1} which also works with the spanning tree ${\mathcal T}_n$ of Section~\ref{redef}. In particular, the switch $\tau_{n,d,z}$ is usable by Lemma~\ref{ult}. The effective shifts of the switches $\tau_{n,1}$, $\tau_{n,2}$ and $\tau_{n,d,z}$ used in the proof, i.e. the signs $\gamma_1$, $\gamma_2$ and $\gamma_d$ in displays (\ref{28a})-(\ref{28c}), can now be computed explicitly. Specifically, from Lemmas~\ref{tai}, \ref{11} and \ref{ult}, it is obtained: \begin{eqnarray}\label{gam}\gamma_1=(-1)(-1)=+1,\hspace*{1cm}\gamma_2=(+1)(+1)=(+1),\hspace*{1cm}\gamma_d=(-1)(+1)=-1.\end{eqnarray} In each product in (\ref{gam}), the first factor is $-1$ if and only if  the switch is $f^{-1}$-conformal and the second factor is $-1$ if and only if  the switch is reversed.

\begin{proof}
The input of the algorithm is the integer $n\ge 1$, the initial combination $x'$ and the desired shift $\bar{s}$ coprime to $2n+1$. Then, $C_n$ mod $2n+1$ is computed in $O(n^2)$ time via Segner's recurrence relation, namely $$C_0=1\;\mbox{   and   }\;C_{n+1}=\sum_{i=0}^{n}C_iC_{n-i}.$$ By (\ref{Cat}), this yields the shift $s$ of the flip sequence obtained from gluing without modifications. A test of whether $s$ is coprime to $2n+1$ proceeds, followed by the computation of one or two switches such that the shift $s'$ of the modified flip sequence $s'$, via (\ref{28a})-(\ref{28c}) and (\ref{gam}), is coprime to $2n+1$. In particular, the definition of $d$ in (\ref{28c}) involves computing the prime factorization of $2n+1$. From $s'$, the {\it scaling factor} $(s')^{-1}\bar{s}$ is computed as well as the corresponding initial combination $x$, such that $x'$ is obtained from $x$ by permuting columns according to the rule $i\mapsto (s')^{-1}\bar{s}i$ (applying the inverse permutation).  These remaining initial steps can be performed in $O(n)$ time. All further computations are then performed with $x$, and whenever a flip position is computed for $x$, it is scaled by $(s')^{-1}\bar{s}$, before applying it to $x'$. 

To decide whether to perform an $f$-step or a pull/push step, the following computations are performed on the current plane tree $T=[t(x)]$, following steps (T1)-)(T3), Section~\ref{redef}:
\begin{enumerate}
\item compute a centroid $c$ of $T$ and its potencial $\phi(c)$, as in (T1), in time $O(n)$ (see~\cite{Kang});
\item compute the lexicographic subtree ordering, as in (T1), in time $O(n)$; if the centroid is unique, this is done via Booth's algorithm~\cite{Boo}; specifically, to compute the lexicographically smallest clockwise (differing from ccw in \cite{mmm}) ordering $(t_1,\ldots,t_k)$ of the $c$-subtrees of $T$, insert $-1$'s as separators between the bitstrings $t_i$; this takes to the string $z:=(-1,t_1,-1,,\ldots,-1,t_k)$; such trick makes Booth's algorithm return a cyclic rotation of $z$ that starts with -1; it follows that such rotation is also the one minimizing the cyclic subtree ordering $(t_1,\ldots,t_k)$;  
\item compute a $c$-subtree of $T$ and one of its leaves, as in steps (T2)-(T3), Section~\ref{redef}.
\end{enumerate} 
Overall, the decision of which type of step to perform next takes time $O(n)$ to compute.
Whenever a switch appears in the course of the algorithm, detectable in time $O(n)$, a modified flip as in (\ref{25a})-(\ref{25c}) is performed. Each time this occurs, the position $\ell(c)$ has to be recomputed, with $T=[t(x)]$ unchanged.

In sum, the algorithm runs in time $O(n)$ in each step, using $O(n)$ memory all the time, and it requires time $O(n^2)$ for initialization.
\end{proof}

\end{document}